\documentclass{article}

\usepackage{amssymb, amsmath, amsthm, array}

\usepackage[left=1.5in, top=1.3in, right=1.3in, bottom=1.2in]{geometry}

\newtheorem{theorem}{Theorem}[section]
\newtheorem{lemma}[theorem]{Lemma}
\newtheorem{proposition}[theorem]{Proposition}
\newtheorem{corollary}[theorem]{Corollary}
\newtheorem{definition}[theorem]{Definition}

\begin{document}

\title{\bf Alternate Definitions of Vector Space Dimension and Module Rank Using Isomorphisms}
\author{Julia Maddox \\
        julia@ou.edu \\
		David and Judi Proctor Department of Mathematics \\
		University of Oklahoma \\
		601 Elm Avenue, Room 423 \\
		Norman, OK 73019-3103}
\date{May 2023}

\maketitle

\abstract{The standard definition of the dimension of a vector space or rank of a module states that dimension or rank is equal to the cardinality of any basis, which requires an understanding of the concepts of basis, generating set, and linear independence. We pose new definitions for the dimension of a vector space, called the isomorphic dimension, and for the rank of a module, called the isomorphic rank, using isomorphisms. In the finite case, for a vector space $V$ over field $F$, its isomorphic dimension is equal $n$ if and only if there exists a linear isomorphism from $F^n$ to $V$. For a module $M$ over the commutative ring $R$ with identity, its isomorphic rank is equal to $n$ if and only if there exists an $R$-module isomorphism from $R^n$ to $M$. There are similar definitions in the infinite cases. These isomorphic definitions do not require the concepts of basis, generating set, and linear independence. This approach allows for some fundamental linear algebra and module theory results to be seen more easily or to be proven more similarly to other algebraic proofs involving isomorphisms and homomorphisms and provides an alternate educational approach to dimension and rank.}

\vspace{0.2in}
classification: 15-01, 15A03, 97H60, 13C05

\vspace{0.2in}
keywords: dimension; rank; basis.

\section{Introduction}

In the history of the epistemology of abstract vector space theory, the definition of dimension has developed to represent an enumeration of a set of elements satisfying certain conditions. Those conditions have evolved from various presentations including a maximal linearly independent set or a minimal generating set to the standard approach employed today where these notions meet. 

The following is a brief recounting of some of the history of basis and dimension. In 1844 in \cite{grassmann1844} and subsequently in 1862 in\cite{grassmann1862}, Grassmann defines dimension in terms of the theory of extension and creates an exchange method from which the exchange theorem can be derived. In 1888 in \cite{peano}, Peano provided an early axiomatic definition of what we now call a vector space. He defined dimension as the maximum number of independent elements in a system and without connecting it to the notion of a generating set. In 1893 in the fourth edition of \cite{dirichlet}, Dedekind contributed a supplement defining a basis in the context of field extensions as an irreducible, or linearly independent, generating set, and he proves a result equivalent to the exchange theorem and presents the size of the basis as an invariant number. In 1909 in \cite{burali-forti marcolongo}, Burali-Forti and Marcolongo continued Peano's work and defined dimension as the maximum number of independent elements as well. In 1910 in \cite{steinitz}, Steinitz gave an explicit definition of linear dependence and defined the dimension of a finite field extension to be the maximum number of linearly independent elements. He considered a basis to be a generating set for which every element has a unique linear expression. In his paper, he proved that the dimension is the minimal size of a generating set, provided results equivalent to the exchange theorem, and included in his theorems the invariance of the size of a basis. In 1918 in \cite{weyl}, Weyl defined linear vector-manifolds in an axiomatic approach, which is similar to the definition of vector spaces, and he similarly defines dimension as the maximum number of linearly independent vectors. In 1930, see \cite{vanderwaerden} and \cite{vanderwaerden2}, van der Waerden published a text with a general approach to modern algebra. In his book, van der Waerden used the exchange theorem and proved results involving bases, the invariance of the number of elements in a basis, and dimension. He also introduced the concepts of linear dependence and basis for modules. His book was based on lectures by Emil Artin and Emmy Noether. In 1921 in \cite{noether}, Emmy Noether defined modules and the basis of a module and considered vector spaces as special cases of modules.
See \cite{dorier1995}, \cite{dorier1996}, and \cite{kleiner} for a more detailed overview of this history.

In the exposition of abstract vector space theory leading up to dimension, linear combinations and spans are first defined. The spanning set and linear independence are subsequently defined. A basis of a vector space is then defined as a linearly independent spanning set of that vector space, or it is equivalently defined as a set of vectors for which all vectors in the vector space can be represented as unique linear combinations of this basis. The size of a linearly independent set is less than or equal to the size of a spanning set, and this is demonstrated to show the number of elements in a basis is invariant. The dimension of a vector space is then the cardinality of every basis of the vector space. From here, it can be proven that any two vector spaces of equal dimension are isomorphic. The linear algebra texts \cite{axler}, \cite{halmos}, \cite{lang}, and \cite{roman} were used for reference, and they are examples of texts that follow the above approach to dimension. The rank of a module is also defined to be the invariant number of elements in a basis, which is a linearly independent generating set for that module. We used texts \cite{berrick keating}, \cite{hungerford}, \cite{lang2}, \cite{dummit foote}, and \cite{roman} as references, and they follow this approach.

In this article, we present a different approach with an alternate definition for vector space dimension called the isomorphic dimension and an alternate definition of for module rank called the isomorphic rank. In the finite case, for a vector space over field $F$, we define the isomorphic dimension of the vector space to be equal $n$ if and only if there exists a linear isomorphism between $F^n$ and that vector space. There is a similar definition in the infinite case. This produces a reordering of the fundamental results of abstract vector space theory involving dimension. While the standard approach relies on sets of vectors of the vector space and elemental arguments, this approach focuses on linear transformations between vector spaces. With this alternate definition, we prove statements with functional arguments. Some arguments still require an appeal to elemental work, but most of these arguments involve vectors as outputs of the standard basis of $F^n$. In the finite case, definition of isomorphic rank of a module over a commutative ring $R$ with identity is equal to $n$ if and only if there exists an $R$-module isomorphism between $R^n$ and that module. There is a similar definition in the infinite case.

The standard definitions of dimension and rank require students to learn the definitions of linear combination, span, spanning set, linear independence, and basis before dimension is mentioned. These definitions are also relatively particular to the context of vector spaces and modules (or similar constructions). These alternate definitions requires students to use general abstract algebra concepts of function, homomorphism (linearity in the case of vector spaces), injective, surjective, bijective, and isomorphism along with some elementary matrix theory.

With the isomorphic definition, $F^n$ becomes a focal vector space and $R^n$ becomes a focal module, for field $F$ and commutative ring with identity $R$.
When defining dimension in terms of isomorphisms with $F^n$ or $R^n$, the isomorphic dimension or rank of the zero vector space or module is equal to zero using the isomorphism $f:F^0=\{{\bf 0}\} \to V=\{{\bf 0}_V\}$ such that $f({\bf 0})={\bf 0}_V$ or $f:R^0=\{{\bf 0}\} \to M=\{{\bf 0}_M\}$ such that $f({\bf 0})={\bf 0}_M$. These definitions provide a constructive way of determining the dimension of the zero vector space or module and do not rely on the abstraction of using the empty set as a basis.

In our approach, some fundamental linear algebra and module theory results are seen more easily or are proven more similarly to other algebraic proofs involving isomorphisms and homomorphisms.
We provide the details of the proofs to show the concepts used. We also include the results specific to spanning sets and linearly independent sets using the framework of surjective and injective linear transformations.

\section{Isomorphic Dimension in the Finite Case}

Let $\{e_1,\ldots,e_n\} \subset F^n$ for field $F$ be the standard basis of $F^n$. For every $i=1,\ldots,n$, $e_i$ has each $j$th coordinate equal to $0$ for $j \neq i$ and the $i$th coordinate equal to $1$. Here the name standard basis is used, but we will only use this as a label without any specific meaning attached to the term basis.

The following theorem is a known elementary theorem that is true by the definition of linear transformation and because for any $(x_1,\ldots,x_n) \in F^n$, $(x_1,\ldots,x_n)=x_1e_1+\cdots+x_ne_n$. 

\begin{theorem}
Let $V$ be a vector space over field $F$. Then for any set or multiset $\{v_1,\ldots,v_n\}$ for $v_1,\ldots,v_n \in V$, $f:F^n \to V$ such that
$$f(x_1,\ldots,x_n)=x_1v_1+\cdots+x_nv_n$$
for all $x_1,\ldots,x_n \in F$ is the unique linear transformation such that
$f(e_i)=v_i$ for all $i=1,\ldots,n$.
\end{theorem}

When $n=0$, $F^0=\{{\bf 0}\}$, and the only linear transformation $f:\{{\bf 0}\} \to V$ is defined by $f({\bf 0})={\bf 0}_V$.

One of the difficulties in proving the dimension of a vector space is well-defined is in showing that any two bases of $V$ contain the same number of vectors. This is made simpler with our alternate definition and using elementary matrix theory.

\begin{lemma} \label{injlemma}
Let $F$ be a field, and let $n$ and $m$ be nonnegative integers. If $f:F^n \to F^m$ is an injective linear transformation, then $n \leq m$.
\end{lemma}

When $n=0$, we mean $F^0=\{{\bf 0}\}$. 

In the following proof, we will interchangeably write vectors of $F^n$ as $n$-tuples and the corresponding $n \times 1$ columns.

\begin{proof}
If $n=0$, then $0 \leq m$ is trivially true. 

If $m=0$, then $f:F^n \to F^0$ is the linear transformation with $f(x_1,\ldots,x_n)={\bf 0}$ for all $(x_1,\ldots,x_n) \in F^n$, and for positive $n$, $f(1,\ldots,1)={\bf 0}=f(0,\ldots,0)$. This contradicts $f$ being injective. So $n=m=0$ in this case. 

Now consider $n,m>0$, and assume $n>m$.
For any $(x_1,\ldots,x_n) \in F^n$, $$f(x_1,\ldots,x_n)=x_1f(e_1)+\cdots+x_nf(e_n)$$ because $f$ is a linear transformation. Therefore, $$f(x_1,\ldots,x_n)=\begin{bmatrix} f(e_1) & \cdots & f(e_n) \end{bmatrix}\begin{bmatrix} x_1 \\ \vdots \\ x_n \end{bmatrix}.$$ Here $A=\begin{bmatrix} f(e_1) & \cdots & f(e_n) \end{bmatrix}$ is an $m \times n$ matrix in $M(m \times n, F)$.

The kernel of $f$ is equal to $\{X \in F^n|AX={\bf 0}_{m \times 1}\}$. From elementary matrix theory, we know that $\ker(f)=\{X \in F^n|{\rm rref}(A)X={\bf 0}_{m \times 1}\}$, and since $n>m$, the number of independent variables in the solution set of ${\rm rref}(A)X={\bf 0}_{m \times 1}$ is at least one. Therefore, ${\bf 0}_{n \times 1}$ and some nonzero vector in $F^n$ both have output ${\bf 0}_{m \times 1}$, and $f$ is not injective. This is a contradiction, so we must have $n \leq m$.
\end{proof}

\begin{proposition} \label{m<=n}
Let $V$ be a vector space over field $F$. If $f_m:F^m \to V$ is an injective linear transformation and $f_n:F^n \to V$ is a linear isomorphism, then $m \leq n$.
\end{proposition}

\begin{proof}
$f_n^{-1} \circ f_m:F^m \to F^n$ is an injective linear transformation, and by Lemma \ref{injlemma}, $m \leq n$.
\end{proof}

\begin{corollary} \label{m=n}
Let $V$ be a vector space over field $F$. If $f_m:F^m \to V$ is a linear isomorphsim and $f_n:F^n \to V$ is a linear isomorphism, then $m=n$.
\end{corollary}

\begin{proof}
By Proposition \ref{m<=n}, $m \leq n$ and $n \leq m$, so $m=n$.
\end{proof}

If there exists a linear isomorphism between $F^n$ and $V$, then $n$ is unique.

\begin{definition}
Let $V$ be a vector space over field $F$. The {\bf isomorphic dimension} of $V$ is $n$, labeled $\dim_{I}(V)=n$, if and only if there exists a linear isomorphism $f:F^n \to V$.
\end{definition}

When $n=0$, $F^0=\{{\bf 0}\}$, and $f:\{{\bf 0}\} \to V$ is a linear isomorphism if and only if $f({\bf 0})={\bf 0}_V$ and $V=\{{\bf 0}_V\}$.
Therefore $\dim_{I}(V)=0$ if and only if $V=\{{\bf 0}_V\}$. 

For any nonnegative integer $n$, $\dim_I(F^n)=n$ using the identity function on $F^n$.

Let $V$ be a vector space over field $F$. If $\dim_{I}(V)=n$ for some nonnegative integer $n$, then $V$ is {\bf finite-dimensional}, and if for any nonnegative integer $n$, there does not exists an isomorphism $f:F^n \to V$, then $V$ is {\bf infinite-dimensional}.

The isomorphic definition of dimension allows for a short proof of the following.

\begin{theorem} \label{isoequaldim}
Let $V$ and $W$ be vector spaces over field $F$, and let $V$ be finite-dimensional. There exists linear isomorphism $f:V \to W$ if and only if $W$ is finite-dimensional and $\dim_{I}(V)=\dim_{I}(W)$.
\end{theorem}

\begin{proof}
Let $\dim_{I}(V)=n$. There exists linear isomorphism $f_V:F^n \to V$. 

If $f:V \to W$ is a linear isomorphism, then $f \circ f_V:F^n \to W$ is a linear isomorphism. Therefore $W$ is finite-dimensional and $\dim_I(W)=n$.

If $W$ is finite-dimensional and $\dim_I(W)=n$, then there exists linear isomorphism $f_W:F^n \to W$, and $f_W \circ f_V^{-1}:V \to W$ is a linear isomorphism.
\end{proof}

We also say $V$ and $W$ are isomorphic if and only if there exists a linear isomorphism between them. If $V$ is finite-dimensional, $V$ and $W$ are isomorphic if and only if $\dim_I(V)=\dim_I(W)$.

For every pair of nonnegative integers $(i,j)$, define $p_i^j:F^i \to F^j$ such that for any $(x_1,\ldots,x_i) \in F^i$, $p_i^j(x_1,\ldots,x_i)=(x_1,\ldots,x_i,0,\ldots,0)$ if $i \leq j$ and $p_i^j(x_1,\ldots,x_i)=(x_1,\ldots,x_j)$ if $i \geq j$.

\begin{lemma} \label{addinjective}
Let $V$ be a vector space over field $F$. Let $f_n:F^n \to V$ be a linear transformation with $v_i=f_n(e_i)$ for all $i=1,\ldots,n$.
If $f_n$ is injective but not surjective, then there exists $v_{n+1} \in V$ with $v_{n+1} \notin {\rm im}(f_n)$, and for any such $v_{n+1} \in V$, there exists an injective linear transformation $f_{n+1}:F^{n+1} \to V$ with $f_{n+1}(e_i)=v_i$ for all $i=1,\ldots,n+1$, $f_{n+1} \circ p_n^{n+1} = f_n$, and ${\rm im}(f_n) \subsetneq {\rm im}(f_{n+1})$.
\end{lemma}

In this statement, $\{e_1,\ldots,e_n\}$ is the standard basis of $F^n$ when used in the context of $f_n:F^n \to V$, and $\{e_1,\ldots,e_{n+1}\}$ is the standard basis of $F^{n+1}$ when used in the context of $f_{n+1}:F^{n+1} \to V$.

\begin{proof}
If $f_n:F^n \to V$ is not surjective, there exists $v_{n+1} \in V$ such that $v_{n+1} \notin {\rm im}(f_n)$. Then for any $a \in F$, $av_{n+1} \in {\rm im}(f_n)$ if and only if $a=0$ because ${\rm im}(f_n)$ is a subspace.

Define $f_{n+1}:F^{n+1} \to V$ to be the linear transformation such that $f_{n+1}(e_i)=v_i$ for all $i=1,\ldots,n+1$. It follows that $f_{n+1} \circ p_n^{n+1} = f_n$ and ${\rm im}(f_n) \subset {\rm im}(f_{n+1})$ with ${\rm im}(f_n) \neq {\rm im}(f_{n+1})$ because $v_{n+1} \notin {\rm im}(f_n)$ but $v_{n+1} \in {\rm im}(f_{n+1})$.

For any $(x_1,\ldots,x_{n+1}) \in \ker(f_{n+1})$, $$f_n(x_1,\ldots,x_n)+x_{n+1}v_{n+1}={\bf 0},$$ which implies $f_n(x_1,\ldots,x_n)=-x_{n+1}v_{n+1}$. Therefore, $f_n(x_1,\ldots,x_n)={\bf 0}_V$ and $x_{n+1}=0$. Since $f_n$ is injective, $(x_1,\ldots,x_n)=(0,\ldots,0)$, and $(x_1,\ldots,x_n,x_{n+1})=(0,\ldots,0,0)$. This proves $f_{n+1}$ is injective.
\end{proof}

\begin{corollary} \label{cor:inj-iso}
Let $V$ be a finite-dimensonal vector space over field $F$ with $\dim_I(V)=n$. If $f:F^n \to V$ is an injective linear transformation, then $f$ is a linear isomorphism.
\end{corollary}

\begin{proof}
If $f$ is not surjective, then by Lemma \ref{addinjective}, there exists injective linear transformation $f_{n+1}:F^{n+1} \to V$. By Proposition \ref{m<=n}, this is a contradiction, so $f$ must be surjective and a linear isomorphism.
\end{proof}

\begin{theorem} \label{injseqthm}
Let $V$ be a vector space over field $F$. 

$V$ is finite-dimensional with dimension $n$ if and only if there exists a sequence of injective linear transformations:
$$f_{0}:F^0 \to V, \, f_{1}:F^1 \to V, \ldots, \, f_{n}:F^n \to V$$
such that $f_{k+1} \circ p_k^{k+1} = f_k$ and ${\rm im}(f_{k}) \subsetneq {\rm im}(f_{k+1})$ for any $k \leq n-1$ and ${\rm im}(f_n)=V$.

$V$ is infinite-dimensional if and only if there exists an infinite sequence of injective linear transformations:
$$f_{0}:F^0 \to V, \, f_{1}:F^1 \to V, \ldots, \, f_{n}:F^n \to V, \ldots$$
such that $f_{k+1} \circ p_k^{k+1} = f_k$ and ${\rm im}(f_{k}) \subsetneq {\rm im}(f_{k+1})$ for any $k \geq 0$.
\end{theorem}

\begin{proof}
For any vector space $V$, there exists the injective linear transformation $f_{0}:F^0 \to V$ defined by $f({\bf 0})={\bf 0}_V$. If $V=\{{\bf 0}_V\}$, then $f_{0}$ satisfies the statement of the theorem. If $V$ is nonzero, there exists nonzero $v_1 \in V$, and there exists the injective linear transformation $f_{1}:F^1 \to V$ defined by $f(ae_1)=av_1$ for any $a \in F$. By Lemma \ref{addinjective}, we can continue building injective linear transformations. The process terminates if we reach an isomorphism.

For any such sequence, using Proposition \ref{m<=n}, $V$ is finite-dimensional with dimension $n$ if and only if this sequence ends with isomorphism $f_n:F^n \to V$, and $V$ is infinite-dimensional if and only if there is no such isomorphism, and the sequence never ends.
\end{proof}

This lemma creates a chain of subspaces of $V$ in the form of the images of the injective linear transformations, and it is reminiscent of the length of a module.

\begin{theorem} \label{subspacedim}
Let $V$ be a finite-dimensional vector space over field $F$, and let $U$ be a subspace of $V$. Then $U$ is finite-dimensional and $\dim_I(U) \leq \dim_I(V)$.
\end{theorem}

\begin{proof}
Let $\dim_I(V)=n$. There exists a linear isomorphism $f:F^n \to V$. Let $h:F^k \to U$ be an injective linear transformation for some nonnegative integer $k$. We know such an injective linear transformation exists for at least $k=0$ defined by $h({\bf 0})={\bf 0}_V$.

Let ${\rm inc}:U \to V$ be the injective linear transformation ${\rm inc}(u)=u$ for any $u \in U$. Then $f^{-1} \circ {\rm inc} \circ h:F^k \to F^n$ is an injective linear transformation. By Lemma \ref{injlemma}, $k \leq n$. Therefore $U$ is finite-dimensional by Theorem \ref{injseqthm}, and $\dim_I(U) \leq \dim_I(V)$.
\end{proof}

\begin{theorem} \label{U=V}
Let $V$ be a finite-dimensional vector space over field $F$, and let $U$ be a subspace of $V$. $U=V$ if and only if $\dim_I(U)=\dim_I(V)$.
\end{theorem}

\begin{proof}
Since $U$ is a subspace of $V$, $U \subset V$.

If $U=V$, then $\dim_I(U)=\dim_I(V)$ trivially.

Let $\dim_I(U)=\dim_I(V)=n$. There exists a linear isomorphism $f:F^n \to U$ because $\dim_I(U)=n$. Let ${\rm inc}:U \to V$ be the injective linear transformation ${\rm inc}(u)=u$ for any $u \in U$. Then ${\rm inc} \circ f:F^n \to V$ is an injective linear transformation. By Corollary \ref{cor:inj-iso}, ${\rm inc} \circ f$ is a linear isomorphism. Since $f$ is also a linear isomorphism, ${\rm inc}=({\rm inc} \circ f) \circ f^{-1}$ is a linear isomorphism, and ${\rm inc}$ is surjective. Thus $U=V$.
\end{proof}

\begin{lemma} \label{injseqUlemma}
Let $V$ be a finite-dimensional vector space over field $F$ with dimension $n$, and let $U$ be a subset of $V$.
$U$ is a subspace of $V$ with dimension $k$ if and only if there exists a sequence of injective linear transformations:
$$f_{0}:F^0 \to V, \, f_{1}:F^1 \to V, \ldots, \, f_{n}:F^n \to V$$
such that $f_{m+1} \circ p_m^{m+1} = f_m$ and ${\rm im}(f_{m}) \subsetneq {\rm im}(f_{m+1})$ for any $m \leq n-1$, ${\rm im}(f_n)=V$, and ${\rm im}(f_k)=U$.
\end{lemma}

\begin{proof}
If ${\rm im}(f_k)=U$ for some nonnegative integer $k \leq n$, then $U$ is a subspace of $V$ because the image of a linear transformation is a subspace of the codomain. As a consequence of the First Isomorphism Theorem $U$ is isomorphic to $F^k$ and $\dim_I(U)=k$.

If $U$ is a subspace of $V$ with $\dim_I(U)=k$, then $k \leq n$. By Theorem \ref{injseqthm}, there exists a sequence of injective linear transformations:
$$h_{0}:F^0 \to U, \, h_{1}:F^1 \to U, \ldots, \, h_{k}:F^n \to U$$
such that $h_{m+1} \circ p_m^{m+1} = h_m$ and ${\rm im}(h_{m}) \subsetneq {\rm im}(h_{m+1})$ for any $m \leq k-1$ and ${\rm im}(h_k)=U$. 

Let ${\rm inc}:U \to V$ be the injective linear transformation such that ${\rm inc}(u)=u$ for any $u \in U$, and define $f_i={\rm inc} \circ h_i:F^i \to V$. Therefore, ${\rm im}(f_k)=U$. Using Lemma \ref{addinjective}, we continue constructing injective linear transformations, until we reach the linear isomorphism $f_n:F^n \to V$. 
\end{proof}

\begin{lemma} \label{isolemma}
Let $V$ and $W$ be vector spaces over field $F$, let $f:V \to W$ be a linear isomorphism, and let $U$ be a subspace of $V$. Then $\bar{f}:V/U \to W/f(U)$ such that $\bar{f}(v+U)=f(v)+f(U)$ for any $v \in V$ is a linear isomorphism.
\end{lemma}

\begin{proof}
$\bar{f}$ is a well-defined linear transformation that is injective and surjective. Thus it is a linear isomorphism.
\end{proof}

\begin{theorem} \label{dim(V/U)}
Let $V$ be a finite-dimensional vector space over field $F$, and let $U$ be a subspace of $V$. $V/U$ is finite-dimensional with $\dim_I(V/U)=\dim_I(V)-\dim_I(U).$
\end{theorem}

\begin{proof}
Let $\dim_I(V)=n$. If $U$ is a subspace of $V$, then $U$ is finite-dimensional with $\dim_I(U)=k \leq n$.
By Lemma \ref{injseqUlemma}, there exists a sequence of injective functions:
$$f_{0}:F^0 \to V, \, f_{1}:F^1 \to V, \ldots, \, f_{n}:F^n \to V$$
such that $f_{m+1} \circ p_m^{m+1} = f_m$ and ${\rm im}(f_{m}) \subsetneq {\rm im}(f_{m+1})$ for any $m \leq n-1$, ${\rm im}(f_n)=V$, and ${\rm im}(f_k)=U$.
Then $f_n(p_k^n(F^k))=U$ for $p_k^n:F^k \to F^n$ such that $p_k^n(x_1,\ldots,x_k)=(x_1,\ldots,x_k,0,\ldots,0)$ for all $(x_1,\ldots,x_k) \in F^k$. 

By Lemma \ref{isolemma}, there exists an isomorphism between $F^n/p_k^n(F^k)$ and $V/f_n(p_k^n(F^k))$ with $f_n(p_k^n(F^k))=U$. By Theorem \ref{isoequaldim}, $\dim_I(V/U)=\dim_I(F^n/p_k^n(F^k))$, which equals $n-k$ using the linear isomorphism $h:F^{n-k} \to F^n/p_k^n(F^k)$, such that $h(x_{k+1},\ldots,x_{n})=(0,\ldots,0,x_{k+1},\ldots,x_n)+p_k^n(F^k)$ for all $(x_{k+1},\ldots,x_n) \in F^{n-k}$.
\end{proof}

\begin{corollary} \label{cor:surj-iso}
Let $V$ be a finite-dimensional vector space over field $F$ with $\dim_I(V)=n$. If $f:F^n \to V$ is a surjective linear transformation, then $f$ is a linear isomorphism.
\end{corollary}

\begin{proof}
By the First Isomorphism Theorem, $F^n/\ker(f)$ is isomorphic to $V$. By Theorem \ref{isoequaldim}, $\dim_I(F^n/\ker(f))=\dim(V)=n$. By Theorem \ref{dim(V/U)}, $\dim_I(F^n/\ker(f))=\dim_I(F^n)-\dim_I(\ker(f))=n-\dim_I(\ker(f))$, so $n-\dim_I(\ker(f))=n$. Therefore $\dim_I(\ker(f))=0$ and $\ker(f)=\{{\bf 0}\}$, so $f$ is injective and a linear isomorphism.
\end{proof}

\begin{theorem} \label{thm:inj:dim(V)<=dim(W)}
Let $V$ and $W$ be vector spaces over field $F$ and such that $W$ is finite-dimensional. There exists an injective linear transformation $f:V \to W$ if and only if $V$ is finite-dimensional and $\dim_I(V) \leq \dim_I(W)$.
\end{theorem}

\begin{proof}
Let $\dim_I(W)=m$. There exists linear isomorphism $f_W:F^m \to W$. 

If $V$ is finite-dimensional, then let $\dim_I(V)=n$ and there exists a linear isomorphism $f_V:F^n \to V$. If $n \leq m$, then $p_n^m:F^n \to F^m$ is the injective linear transformation $p_n^m(x_1,\ldots,x_n)=(x_1,\ldots,x_n,0,\ldots,0)$ for any $(x_1,\ldots,x_n) \in F^n$, and $f_W \circ p_n^m \circ f_V^{-1}:V \to W$ is an injective linear transformation.

If $f:V \to W$ is an injective linear transformation, then by the First Isomorphism Theorem, $V$ is isomorphic to ${\rm im}(f)$. Since ${\rm im}(f)$ is a subspace of finite-dimensional $W$, it is finite-dimensional and its dimension is less than or equal to $m$ by Theorem \ref{subspacedim}. By Theorem \ref{isoequaldim}, $V$ is finite-dimensional with $\dim_I(V)=\dim_I({\rm im}(f))$. Therefore $\dim_I(V) \leq \dim_I(W)$.
\end{proof}

\begin{theorem} \label{thm:surj:dim(V)>=dim(W)}
Let $V$ and $W$ be vector spaces over field $F$ and such that $V$ be finite-dimensional. There exists surjective linear transformation $f:V \to W$ if and only if $W$ is finite-dimensional and $\dim_I(V) \geq \dim_I(W)$.
\end{theorem}

\begin{proof}
Let $\dim_I(V)=n$. There exists linear isomorphism $f_V:F^n \to W$.

If $W$ is finite-dimensional, then let $\dim_I(W)=m$ and there exists linear isomorphism $f_W:F^m \to W$. If $n \geq m$, then $p_n^m:F^n \to F^m$ is the surjective linear transformation $p_n^m(x_1,\ldots,x_n)=(x_1,\ldots,x_m)$ for all $(x_1,\ldots,x_n) \in F^n$, and $f_W \circ p_n^m \circ f_V^{-1}:V \to W$ is a surjective linear transformation.

If $f:V \to W$ is a surjective linear transformation, then by the First Isomorphism Theorem, $W$ is isomorphic to $V/\ker(f)$. Since $\ker(f)$ is a subspace of finite-dimensional $V$, $V/\ker(f)$ is finite-dimensional with $\dim_I(V/\ker(f))=\dim_I(V)-\dim_I(\ker(f))$ by Theorem \ref{dim(V/U)}. By Theorem \ref{isoequaldim}, $W$ is finite-dimensional with $\dim_I(W)=\dim_I(V)-\dim_I(\ker(f))$. Therefore $\dim_I(V) \geq \dim_I(W)$.
\end{proof}

The following theorem is the Rank Nullity Theorem using the definition of isomorphic dimension.

\begin{theorem}
Let $V$ and $W$ be vector spaces over field $F$, let $V$ be finite-dimensional, and let $f:V \to W$ be a linear transformation. Then $\ker(f)$ and ${\rm im}(f)$ are finite-dimensional and

$$\dim_I(\ker(f))+\dim_I({\rm im}(f))=\dim_I(V).$$
\end{theorem}

\begin{proof}
Since $\ker(f)$ is a subspace of $V$ and $V$ is finite-dimensional, then $\ker(f)$ is finite-dimensional by Theorem \ref{subspacedim}.
By Theorem \ref{dim(V/U)}, $V/\ker(f)$ is finite-dimensional and $\dim_I(V/\ker(f))=\dim_I(V)-\dim_I(\ker(f))$.
By the First Isomorphism Theorem, $\bar{f}:V/\ker(f) \to {\rm im}(f)$ is a linear isomorphism, and by Theorem \ref{isoequaldim}, ${\rm im}(f)$ is finite-dimensional and $\dim_I(V/\ker(f))=\dim_I({\rm im}(f))$. Therefore $\dim_I(V)-\dim_I(\ker(f))=\dim_I({\rm im}(f))$.
\end{proof}

\section{Isomorphic Basis of a Vector Space in the Finite Dimension Case}

In a similar fashion to our latest definition of dimension, we can also construct a new definition for basis that is equivalent to the algebraic definition of a basis.

\begin{definition}
Let $V$ be a finite-dimensional vector space over field $F$. The set (or multiset) $\{v_1,\ldots,v_n\}$ for $v_1,\ldots,v_n \in V$ is an {\bf isomorphic basis} of $V$ if and only if the linear transformation $f:F^n \to V$ such that $f(e_i)=v_i$ for all $i=1,\ldots,n$ is an isomorphism.
\end{definition}

The standard basis of $F^n$, $\{e_1,\ldots,e_n\}$, is also an isomorphic basis of $F^n$ according to this definition using the identity map.

Let $V$ be a vector space over field $F$. Since a linear transformation with domain $F^n$ is uniquely determined by the outputs of the standard basis, for $\dim_I(V)=n$, there is a one-to-one correspondence between the linear isomorphisms $f:F^n \to V$ and the ordered bases of $V$.

\begin{theorem}
Let $V$ be a vector space over field $F$. Then $V$ is finite-dimensional with $\dim_I(V)=n$ if and only if there exists an isomorphic basis of $V$ containing exactly $\dim_I(V)$ vectors.
\end{theorem}

\begin{proof}
If $\dim_I(V)=n$, then there exists a linear isomorphism $f:F^n \to V$. Then $\{f(e_1),\ldots,f(e_n)\}$ is an isomorphic basis of $V$ by definition.

If $\{v_1,\ldots,v_n\}$ is an isomorphic basis of $V$, then the linear transformation $f:F^n \to V$ defined by $f(e_i)=v_i$ for all $i=1,\ldots,n$ is an isomorphism. Since there exists a linear isomorphism between $F^n$ and $V$, $V$ is finite-dimensional and $\dim_I(V)=n$.
\end{proof}

If $V$ is a finite-dimensional vector space, $\dim_I(V)$ is unique, so we get the following corollary.

\begin{corollary}
Let $V$ be a finite-dimensional vector space over field $F$. Every isomorphic basis of $V$ contains exactly $\dim_I(V)$ vectors.
\end{corollary}

When $\dim_I(V)=0$, the only isomorphic basis of $V$ is the empty set.

\begin{definition}
Let $V$ be a vector space over field $F$. For the set (or multiset) $\{v_1,\ldots,v_n\}$ with $v_1,\ldots,v_n \in V$, a {\bf linear combination} of $\{v_1,\ldots,v_n\}$ is an output of the linear transformation $f:F^n \to V$ such that $f(e_i)=v_i$ for all $i=1,\ldots,n$. 
\end{definition}

This means that a linear combination is $f(x_1,\ldots,x_n)=x_1v_1+\cdots+x_nv_n$ for some $x_1,\ldots,x_n \in F$, which is the same as the usual definition.

\begin{theorem}
Let $V$ be a finite-dimensional vector space over field $F$, and let $\mathcal{B}=\{v_1,\ldots,v_n\}$ for some $v_1,\ldots,v_n \in V$. Then $\mathcal{B}$ is an isomorphic basis of $V$ if and only if every vector in $V$ is a unique linear combination of $\mathcal{B}$.
\end{theorem}

\begin{proof}
Let $f:F^n \to V$ be the linear transformation such that $f(e_i)=v_i$ for all $i=1,\ldots,n$. 

If $\mathcal{B}$ is an isomorpic basis, then $f$ is an isomorphism and every vector in $V$ is a unique output for some input $(x_1,\ldots,x_n) \in F^n$, which means every vector in $V$ is a unique linear combination $x_1v_1+\cdots+x_nv_n$ for that $(x_1,\ldots,x_n) \in F^n$. 

If every vector in $V$ is a unique linear combination of $\mathcal{B}$, then $f$ is surjective because every vector in $V$ is an output of $f$, and $f$ is injective because $f(0,\ldots,0)=0v_1+\cdots+0v_n={\bf 0}_V$ is unique, which implies $\ker(f)=\{{\bf 0}\}$. Thus $f$ is an isomorphism, and $\mathcal{B}$ is an isomorphic basis.
\end{proof}

\section{Surjective and Injective Sets in the Finite Dimension Case}

\begin{definition}
Let $V$ be a vector space over field $F$. For the set (or multiset) $S=\{v_1,\ldots,v_n\}$ for $v_1,\ldots,v_n \in V$, the {\bf span} of $S$ is the image of the linear transformation $f:F^n \to V$ such that $f(e_i)=v_i$ for all $i=1,\ldots,n$, or in other words, ${\rm span}(S)={\rm im}(f)$. 
\end{definition}

This means that ${\rm span}(v_1,\ldots,v_n)=\{x_1v_1+\cdots+x_nv_n|x_1,\ldots,x_n \in F\}$, which is the same as the usual definition.

\begin{definition}
Let $V$ be a vector space over field $F$. The set (or multiset) $\{v_1,\ldots,v_n\}$ for $v_1,\ldots,v_n \in V$ is a {\bf surjective set} in $V$ if and only if the linear transformation $f:F^n \to V$ such that $f(e_i)=v_i$ for all $i=1,\ldots,n$ is surjective.
\end{definition}

Therefore, $S=\{v_1,\ldots,v_n\}$ is a surjective set in $V$ if and only if ${\rm span}(S)=V$. The linear transformation $f:F^n \to V$ defined by $f(e_i)=v_i$ for all $i=1,\ldots,n$ is also surjective if and only if $S$ is a spanning set of $V$. Therefore the definitions of surjective set in $V$ and spanning set of $V$ are equivalent (for finite-dimensional $V$ as currently defined). 

\begin{theorem}
Let $V$ be a vector space over field $F$, and let $S=\{v_1,\ldots,v_n\}$ be a set or multiset of vectors in $V$. Then $S$ is a spanning set of $V$ if and only if $S$ is a surjective set in $V$.
\end{theorem}

We can prove the results involving spanning sets for surjective sets using the definition of surjective sets.

\begin{theorem}
Let $V$ be a vector space over field $F$. $V$ is finite-dimensional if and only if there exists a finite surjective set in $V$.
\end{theorem}

\begin{proof}
If $V$ is finite-dimensional with $\dim_I(V)=n$, then there exists a linear isomorphism $f:F^n \to V$, and $\{f(e_1),\ldots,f(e_n)\}$ is a surjective set in $V$.

If $\{v_1,\ldots,v_n\}$ is a surjective set in $V$, then there exists a surjective linear transformation $g:F^n \to V$ such that $g(e_i)=v_i$ for all $i=1,\ldots,n$. By the First Isomorphism Theorem, there exists a linear isomorphism between $F^n/\ker(g)$ and $V$. Since $F^n/\ker(g)$ is finite-dimensional by Theorem \ref{dim(V/U)}, $V$ is finite-dimensional by Theorem \ref{isoequaldim}.
\end{proof}

\begin{theorem} \label{surjset>=dim(V)}
Let $V$ be a finite-dimensional vector space over field $F$, and let $\{v_1,\ldots,v_m\}$ be a surjective set in $V$. Then $m \geq \dim_I(V)$.
\end{theorem}

\begin{proof}
Let $\dim_I(V)=n$. There exists a linear isomorphism $f:F^n \to V$. Since $\{v_1,\ldots,v_m\}$ is a surjective set in $V$, there exists a surjective linear transformation $g:F^m \to V$ such that $g(e_i)=v_i$ for all $i=1,\ldots,m$.

$f^{-1} \circ g:F^m \to F^n$ is a surjective linear transformation. By the First Isomorphism Theorem, there exists a linear isomorphism between $F^m/\ker(f^{-1} \circ g)$ and $F^n$. Since $F^m/\ker(f^{-1} \circ g)$ is finite-dimensional with dimension $m-\dim_I(\ker(f^{-1} \circ g)$ by Theorem \ref{dim(V/U)}, $m-\dim_I(\ker(f^{-1} \circ g)=n$ by Theorem \ref{isoequaldim} and $m \geq n$.
\end{proof}

\begin{theorem}
Let $V$ be a finite-dimensional vector space over field $F$ with $\dim_I(V)=n$, and let $\{v_1,\ldots,v_n\}$ be a surjective set in $V$, then $\{v_1,\ldots,v_n\}$ is an isomorphic basis of $V$.
\end{theorem}

\begin{proof}
Since $\{v_1,\ldots,v_n\}$ is a surjective set in $V$, the linear transformation $f:F^n \to V$ such that $f(e_i)=v_i$ for all $i=1,\ldots,n$ is surjective. By Corollary \ref{cor:surj-iso}, $f$ is a linear isomorphism, and thus $\{v_1,\ldots,v_n\}$ is an isomorphic basis.
\end{proof}

\begin{theorem}
Let $V$ be a finite-dimensional vector space over field $F$ with $\dim_I(V)=n$, and let $\{v_1,\ldots,v_m\}$ be a surjective set in $V$, then there exists some subset $\{v_{i_1},\ldots,v_{i_n}\} \subset \{v_1,\ldots,v_m\}$ such that $\{v_{i_1},\ldots,v_{i_n}\}$ is an isomorphic basis of $V$.
\end{theorem}

\begin{proof}
Since $\{v_1,\ldots,v_m\}$ is a surjective set in $V$, the linear transformation $f:F^m \to V$ defined by $f(e_i)=v_i$ for all $i=1,\ldots,m$ is surjective. By Theorem \ref{surjset>=dim(V)}, $m \geq n$.

If $\ker(f)=\{{\bf 0}\}$, then $f$ is bijective and an isomorphism, and $m=n$ and $\{v_1,\ldots,v_m\}$ is an isomorphic basis. If $\ker(f) \neq \{{\bf 0}\}$, then $m>n$ and there exists nonzero $(x_1,\ldots,x_m) \in \ker(f)$ with some nonzero coordinate. Without loss of generality, assume $x_m \neq 0$. Define $f_{m-1}:F^{m-1} \to V$ by $f_{m-1}(e_i)=v_i$ for all $i=1,\ldots,m-1$. Since $f$ is surjective, for any $v \in V$, there exists $(a_1,\ldots,a_m) \in F^m$ such that $f(a_1,\ldots,a_m)=a_1v_1+\cdots+a_mv_m=v$. If $a_m=0$, then $f_{m-1}(a_1,\ldots,a_{m-1})=a_1v_1+\cdots+a_{m-1}v_{m-1}=v$. If $a_m \neq 0$, define $(b_1,\ldots,b_m)$ by $$(b_1,\ldots,b_m)=(a_1,\ldots,a_m)-\frac{a_m}{x_m}(x_1,\ldots,x_m)=(a_1-\frac{a_m}{x_m},\ldots,a_{m-1}-\frac{a_m}{x_m},0).$$ Thus $f(b_1,\ldots,b_m)=f(a_1,\ldots,a_m)$ because $(x_1,\ldots,x_m) \in \ker(f)$, and since $b_m=0$, $f_{m-1}(b_1,\ldots,b_{m-1})=v$. Therefore $f_{m-1}$ is surjective.

Continue with this process a total of $m-n$ times until $f_{n}:F^n \to V$ is the linear transformation $f_n(e_i)=v_i$ for all $i=i_1,\ldots,i_n$ and $f_n$ is surjective. By Corollary \ref{cor:surj-iso}, $f_n$ is a linear isomorphism, and $\{v_{i_1},\ldots,v_{i_n}\}$ is a basis.
\end{proof}

\begin{definition}
Let $V$ be a vector space over field $F$. The set (or multiset) $\{v_1,\ldots,v_n\}$ for $v_1,\ldots,v_n \in V$ is an {\bf injective set} in $V$ if and only if the linear transformation $f:F^n \to V$ such that $f(e_i)=v_i$ for all $i=1,\ldots,n$ is injective.
\end{definition}

If $n=0$ and $\{v_1,\ldots,v_n\}=\emptyset$, then $F^0=\{{\bf 0}\}$ and the linear transformation $f:\{{\bf 0}\} \to V$ is defined by $f({\bf 0})={\bf 0}_V$, so $f$ is injective. This implies the empty set is an injective set in $V$.

\begin{theorem}
Let $V$ be a vector space over a field $F$, and let $S=\{v_1,\ldots,v_n\}$ be a set or multiset of vectors in $V$. Then $S$ is a linearly independent set in $V$ if and only if $S$ is an injective set in $V$.
\end{theorem}

\begin{proof}
Let $f:F^n \to V$ be the linear transformation defined by $f(e_i)=v_i$ for all $i=1,\ldots,n$. Then $f(x_1,\ldots,x_n)=x_1v_1+\cdots+x_nv_n$ for all $x_1,\ldots,x_n \in F$ and $(x_1,\ldots,x_n) \in \ker(f)$ if and only if $x_1v_1+\cdots+x_nv_n={\bf 0}_V$. 

If $S$ is a linearly independent set in $V$, then $x_1v_1+\cdots+x_nv_n={\bf 0}_V$ implies $(x_1,\ldots,x_n)=(0,\ldots,0)$, $\ker(f)=\{(0,\ldots,0)\}$, and $f$ is injective, so $S$ is an injective set in $V$. If $S$ is an injective set in $V$, then $f$ is injective and $\ker(f)=\{(0,\ldots,0)\}$, and $x_1v_1+\cdots+x_nv_n={\bf 0}_V$ implies $(x_1,\ldots,x_n)=(0,\ldots,0)$.
\end{proof}

We can prove the results involving linearly independent set for injective sets using the definition of injective sets.

\begin{theorem} \label{injset<=dim(V)}
Let $V$ be a finite-dimensional vector space over field $F$. If $\{v_1,\ldots,v_k\}$ is an injective set in $V$, then $k \leq \dim_I(V)$.
\end{theorem}

\begin{proof}
Let $\dim_I(V)=n$. There exists a linear isomorphism $f:F^n \to V$. Since $\{v_1,\ldots,v_k\}$ is an injective set in $V$, there exists an injective linear transformation $g:F^k \to V$ such that $g(e_i)=v_i$ for all $i=1,\ldots,m$.

$f^{-1} \circ g:F^k \to F^n$ is an injective linear transformation. By Lemma \ref{injlemma}, $k \leq n$.
\end{proof}

\begin{corollary}
Let $V$ be a finite-dimensional vector space over field $F$. If $\{v_1,\ldots,v_k\}$ is an injective set in $V$ and $\{w_1,\ldots,w_m\}$ is a surjective set in $V$, then $k \leq m$.
\end{corollary}

\begin{proof}
By Theorem \ref{surjset>=dim(V)}, $m \geq \dim_I(V)$. By Theorem \ref{injset<=dim(V)}, $k \leq \dim_I(V)$. Therefore $k \leq m$.
\end{proof}

\begin{theorem}
Let $V$ be a finite-dimensional vector space over field $F$ with $\dim_I(V)=n$, and let $\{v_1,\ldots,v_n\}$ be an injective set in $V$, then $\{v_1,\ldots,v_n\}$ is an isomorphic basis of $V$.
\end{theorem}

\begin{proof}
Since $\{v_1,\ldots,v_n\}$ is an injective set in $V$, the linear transformation $f:F^n \to V$ such that $f(e_i)=v_i$ for all $i=1,\ldots,n$ is injective. By Corollary \ref{cor:inj-iso}, $f$ is a linear isomorphism, and $\{v_1,\ldots,v_n\}$ is an isomorphic basis of $V$.
\end{proof}

\begin{theorem}
Let $V$ be a finite-dimensional vector space over field $F$ with $\dim_I(V)=n$, and let $\{v_1,\ldots,v_m\}$ be an injective set in $V$. There exist $v_{m+1},\ldots,v_n \in V$ such that $\{v_1,\ldots,v_n\}$ is an isomorphic basis of $V$.
\end{theorem}

\begin{proof}
Since $\{v_1,\ldots,v_m\}$ is an injective set in $V$, the linear transformation $f:F^m \to V$ such that $f(e_i)=v_i$ for all $i=1,\ldots,m$ is injective.

If $m=n$, then $\{v_1,\ldots,v_m\}$ is an isomorphic basis of $V$ by the previous theorem. If $m<n$, then $f$ is not an isomorphism, and there exists $v_{m+1} \in V$ such that $v_{m+1} \notin {\rm im}(f)$. Define the linear transformation $f_{m+1}:F^{m+1} \to V$ by $f_{m+1}(e_i)=v_i$ for all $i=1,\ldots,m+1$. By Lemma \ref{addinjective}, $f_{m+1}$ is injective. If $m+1=n$, then $\{v_1,\ldots,v_m,v_{m+1}\}$ is an isomorphic basis of $V$ by the previous theorem. If $m+1<n$, continue with this process $n-m$ times. Then $\{v_1,\ldots,v_{m+n-m}=v_n\}$ is an isomorphic basis of $V$ by the previous theorem.
\end{proof}

\section{Isomorphic Rank in the Finite Case}

Let $R$ be a commutative ring with identity, and let $\{e_1,\ldots,e_n\} \subset R^n$ be the standard basis of $R^n$. For every $i=1,\ldots,n$, $e_i$ has each $j$th coordinate equal to $0$ for $j \neq i$ and the $i$th coordinate equal to $1$. Here the name standard basis is used, but we will only use this as a label without any specific meaning attached to the term basis.

The following theorem is a known elementary theorem that is true by the definition of a module homomorphism and because for any $(x_1,\ldots,x_n) \in R^n$, $(x_1,\ldots,x_n)=x_1e_1+\cdots+x_ne_n$. 

\begin{theorem}
Let $R$ be a commutative ring with identity, and let $M$ be a module over field $F$. Then for any set or multiset $\{m_1,\ldots,m_n\}$ for $m_1,\ldots,m_n \in M$, $f:R^n \to M$ such that
$$f(x_1,\ldots,x_n)=x_1m_1+\cdots+x_nm_n$$
for all $x_1,\ldots,x_n \in R$ is the unique module homomorphism such that
$f(e_i)=m_i$ for all $i=1,\ldots,n$.
\end{theorem}

When $n=0$, $R^0=\{{\bf 0}\}$, and the only linear transformation $f:\{{\bf 0}\} \to M$ is defined by $f({\bf 0})={\bf 0}_M$.

\begin{theorem} \label{thmhatf}
Let $R$ be a commutative ring with identity, let $I$ be an ideal of $R$, and let $M$ be an $R$-module. For any surjective $R$-module homomorphism $f:R^n \to M$, $$\hat{f}:R^n/I^n \to M/f(I^n)$$ such that $$\hat{f}((r_1,\ldots,r_n)+I^n)=f(r_1,\ldots,r_n)+f(I^n)$$ for all $(r_1,\ldots,r_n) \in R^n$ is a surjective $R/I$-module homomorphism, and if $f$ is bijective, then $\hat{f}$ is bijective.
\end{theorem}

If $n=0$, all modules are zero modules, and all homomorphisms are defined by sending ${\bf 0}$ to ${\bf 0}$ for the appropriate modules.

\begin{proof}
$I^n$ is an $R$-submodule of $R^n$ because $I$ is an ideal of $R$, and $f(I^n)$ is an $R$-submodule of $M$ because it is the image of an $R$-module homomorphism of an $R$-submodule.

$R^n/I^n$ is an $R/I$-module with scalar multiplication defined as $$(a+I) \cdot ((r_1,\ldots,r_n)+I^n) = a(r_1,\ldots,r_n)+I^n.$$ Scalar multiplication is well-defined because $I$ is an ideal of $R$.

$M/f(I^n)$ is an $R/I$-module with scalar multiplication defined as $$(a+I) \cdot (m+f(I^n)) = am+f(I^n).$$ Scalar multiplication is well-defined because $I$ is an ideal of $R$ and $f$ is a surjective $R$-module homomorphism.

Since $f$ is an $R$-module homomorphism, it follows directly that $\hat{f}$ is an $R$-module homomorphism and an $R/I$-module homomorphism.

Since $f$ is surjective, then for any $m+f(I^n) \in M/f(I^n)$ with $m \in M$, there exists $(r_1,\ldots,r_n) \in R^n$ such that $f(r_1,\ldots,r_n)=m$. Therefore $\hat{f}((r_1,\ldots,r_n)+I^n)=m+f(I^n)$, and $\hat{f}$ is surjective.

Let $f$ be injective. If $(r_1,\ldots,r_n)+I^n \in \ker(\hat{f})$, then $f(r_1,\ldots,r_n)+f(I^n)={\bf 0}_M+f(I^n)$ and $f(r_1,\ldots,r_n) \in f(I^n)$. Since $f$ is injective, $(r_1,\ldots,r_n) \in I^n$ and $\ker(\hat{f})=\{I^n\}$. Therefore $\hat{f}$ is injective.

If $f$ is bijective, then $\hat{f}$ is bijective by the previous two arguments.
\end{proof}

For any $R$-module homomorphism $f$, the definition of $\hat{f}$ in the previous theorem is an $R$-module homomorphism, but we are particularly interested in the case where $f$ is surjective and $\hat{f}$ is an $R/I$-module homomorphism.

For surjective $R$-module homomorphism $f:R^n \to M$, if $f|_{I^n}$ is injective and $\hat{f}$ is an isomorphism, then $f$ is an isomorphism. 

If $I=\{{\bf 0}\}$, then $I^n=\{{\bf 0}\}$ and $f$ and $\hat{f}$ are equivalent. If $R$ is a field, then the only proper ideal of $R$ is $I=\{{\bf 0}\}$, and $f$ and $\hat{f}$ are equivalent linear transformations.

In the following corollary, in which $R$ is a commutative ring with identity, $I$ is an ideal of $R$, and $M$ is an $R$-module, we use $$IM = \{i_1m_1+\cdots+i_km_k|k \in \mathbb{Z}_{\geq 0}, \, i_1,\ldots,i_k \in I, \, m_1,\ldots,m_k \in M\},$$
the set of all finite $I$-linear combinations of elements of module $M$. $IM$ is a submodule of $M$.

\begin{corollary} \label{corhatf}
Let $R$ be a commutative ring with identity, let $I$ be an ideal of $R$, and let $M$ be an $R$-module. For any surjective $R$-module homomorphism $f:R^n \to M$, $$\hat{f}:(R/I)^n \to M/IM$$ such that $$\hat{f}(r_1+I,\ldots,r_n+I)=f(r_1,\ldots,r_n)+IM$$ for all $(r_1,\ldots,r_n) \in R^n$ is a surjective $R/I$-module homomorphism, and if $f$ is bijective, then $\hat{f}$ is bijective.
\end{corollary}

If $n=0$, all modules are zero modules, and all homomorphisms are defined by sending ${\bf 0}$ to ${\bf 0}$ for the appropriate modules.

\begin{proof}
Since $I$ is an ideal of $R$ and $f$ is a surjective $R$-homomorphism, $f(I^n)=IM$.

$(R/I)^n$ and $R^n/I^n$ are isomorphic $R$-modules and isomorphic $R/I$-modules using the following module isomorphism: $(r_1+I,\ldots,r_n+I) \mapsto (r_1,\ldots,r_n)+I^n$. Compose this isomorphism with $\hat{f}$ from the previous theorem and replace $f(I^n)$ with $IM$.
\end{proof}

\begin{theorem} \label{iso:m=n}
Let $R$ be a commutative ring with identity, and let $M$ be an $R$-module. If there exist an $R$-module isomorphism $f:R^n \to M$ and an $R$-module isomorphism $g:R^k \to M$, then $n=k$.
\end{theorem}

\begin{proof}
There exists a maximal ideal $I$ in $R$, and $R/I$ is a field. By the previous corollary, there exist $R/I$-module isomorphisms $\hat{f}:(R/I)^n \to M/IM$ and $\hat{g}:(R/I)^k \to M/IM$. Since $R/I$ is a field, $\hat{f}$ and $\hat{g}$ are linear isomorphisms between vector spaces. Two vector spaces are isomorphic if and only if they have equal dimensions, so $\dim(R/I)^n=\dim(M/IM)=\dim(R/I)^k$ or $n=k$.
\end{proof}

\begin{definition}
Let $R$ be a commutative ring with identity, and let $M$ be an $R$-module. $M$ is {\bf free} and with finite {\bf isomorphic rank} ${\rm rank}_I(M)=n$ if and only if there exists an $R$-module isomorphism $f:R^n \to M$.
\end{definition}

$M=\{{\bf 0}_M\}$ if and only if the $R$-module homomorphism $f:R^0=\{{\bf 0}\} \to M$, such that $f({\bf 0})={\bf 0}_M$, is an isomorphism. Therefore, $M=\{{\bf 0}_M\}$ if and only if ${\rm rank}_I(M)=0$ according to this definition.

For any nonnegative integer $n$, $R^n$ is free with ${\rm rank}_I(R^n)=n$ using the identity map on $R^n$.

\begin{theorem} \label{thmisorank(M)=rank(N)}
Let $R$ be a commutative ring with identity, and let $M$ and $N$ be $R$-modules such that $M$ is free with finite isomorphic rank. Then $M$ and $N$ are isomorphic if and only if $N$ is free with finite isomorphic rank and ${\rm rank}_I(M)={\rm rank}_I(N)$.
\end{theorem}

\begin{proof}
Let ${\rm rank}_I(M)=k$. Then there exists an $R$-module isomorphism $f_M:R^k \to M$.

If $M$ and $N$ are isomorphic, then there exists an $R$-module isomorphism $g:M \to N$, and $g \circ f_M:R^k \to N$ is an $R$-module isomorphism. Therefore $N$ is free with ${\rm rank}_I(N)=k={\rm rank}_I(M)$.

If $N$ is free and ${\rm rank}_I(M)={\rm rank}_I(N)$, then there exists an $R$-module isomorphism $f_N:R^k \to N$, and $f_N \circ f_M^{-1}:M \to N$ is an $R$-module isomorphism. Therefore $M$ and $N$ are isomorphic.
\end{proof}

\begin{theorem} \label{thm:surj:rank(M)>=rank(N)}
Let $R$ be a commutative ring with identity and let $M$ be a free $R$-module such that ${\rm rank}_I(M)=m$, and let $N$ be a free $R$-module such that ${\rm rank}_I(N)=n$. There exists a surjective $R$-module homomorphism $f:M \to N$ if and only if $m \geq n$.
\end{theorem}

\begin{proof}
Since ${\rm rank}_I(M)=m$ and ${\rm rank}_I(N)=n$, there exist an $R$-module isomorphism $g_M:R^m \to M$ and an $R$-module isomorphism $g_N:R^n \to N$.

If there exists a surjective $R$-module homomorphism $f:M \to N$, then $h=g_N^{-1} \circ f \circ g_M:R^m \to R^n$ is a surjective $R$-module homomorphism.

Let $I$ be a maximal ideal of $R$, so that $R/I$ is a field. Then by Corollary \ref{corhatf},
there exists a surjective linear transformation $\hat{h}:(R/I)^m \to R^n/IR^n$. Since $I$ is an ideal of $R$, $IR^n=I^n$. The composition of $\hat{h}$ with the $R/I$-linear isomorphism from $R^n/I^n$ to $(R/I)^n$ defined by $(r_1,\ldots,r_n)+I^n \mapsto (r_1+I,\ldots,r_n+I)$ produces a surjective linear transformation $(R/I)^m \to (R/I)^n$ over field $R/I$. Therefore $\dim((R/I)^m) \geq \dim((R/I)^n)$ or $m \geq n$.

If $m \geq n$, then there exists a surjective $R$-module homomorphism $p:R^m \to R^n$ such that $p(x_1,\ldots,x_m)=(x_1,\ldots,x_n)$, and $f=g_N \circ p \circ g_M^{-1}:M \to N$ is a surjective $R$-module homomorphism.
\end{proof}

\begin{theorem} \label{thm:rank(M/N)}
Let $R$ be a commutative ring with identity, let $M$ be an $R$-module, and let $N$ be a submodule of $M$. If $N$ and $M/N$ are free with finite isomorphic rank, then $M$ is free with finite isomorphic rank and
$${\rm rank}_I(M) = {\rm rank}_I(N)+{\rm rank}_I(M/N).$$
\end{theorem}

\begin{proof}
Let ${\rm rank}_I(N)=n$ and ${\rm rank}_I(M/N)=k$. There exist $R$-module isomorphisms $f_N:R^n \to N$ and $f_{M/N}:R^k \to M/N$. Define $g:M/N \to M$ to be an injective function that sends each element of $M/N$ to a particular coset representative. Define $f_M:R^{n+k} \to M$ as the $R$-module homomorphism such that $f_M(e_i)=f_N(e_i)$ for $i=1,\ldots,n$ and $f_M(e_i)=g(f_{M/N}(e_{i-n}))$ for $i=n+1,\ldots,n+k$. 

Every element $m \in M$ is a sum of any coset representative from $m+N$ and some element of $N$, and since $f_N$ and $f_{M/N}$ are surjective, $f_M$ is surjective. 

For any $(a_1,\ldots,a_{n+k}) \in R^n$, $$f_M(a_1,\ldots,a_{n+k})=a_1f_N(e_1)+\cdots+a_nf_N(e_n)+a_{n+1}g(f_{M/N}(e_1))+\cdots+a_{n+k}g(f_{M/N}(e_k)).$$ For any $(a_1,\ldots,a_{n+k}) \in \ker(f_M)$,
$$a_1f_N(e_1)+\cdots+a_nf_N(e_n)+a_{n+1}g(f_{M/N}(e_1))+\cdots+a_{n+k}g(f_{M/N}(e_k))={\bf 0}_M,$$ which implies $a_{n+1}f_{M/N}(e_1)+\cdots+a_{n+k}f_{M/N}(e_k)={\bf 0}_{M/N}$ in $M/N$. Since $f_{M/N}$ is injective, $a_{n+1}=\cdots=a_{n+k}=0$, and $a_1f_N(e_1)+\cdots+a_nf_N(e_n)={\bf 0}_M$, which implies $a_1=\cdots=a_n=0$ because $f_N$ is injective.

Since $f_M$ is bijective, it is an isomorphism, and $M$ is free with ${\rm rank}_I(M)={\rm rank}_I(N)+{\rm rank}_I(M/N)$.
\end{proof}

\begin{theorem}
Let $R$ be a commutative ring with identity, let $M$ and $N$ be $R$-modules, and let $f:M \to N$ be an $R$-module homomorphism. If $\ker(f)$ and ${\rm im}(f)$ are free with finite isomorphic rank, then $M$ is free with finite isomorphic rank and
$${\rm rank}_I(M) = {\rm rank}_I(\ker(f))+{\rm rank}_I({\rm im}(f)).$$
\end{theorem}

\begin{proof}
By the First Isomorphism Theorem, $M/\ker(f)$ is isomorphic to ${\rm im}(f)$, and by Theorem \ref{thmisorank(M)=rank(N)}, $M/\ker(f)$ is free with finite rank and ${\rm rank}_I(M/\ker(f))={\rm rank}_I({\rm im}(f))$.
By Theorem \ref{thm:rank(M/N)}, $${\rm rank}_I(M)={\rm rank}_I(\ker(f))+{\rm rank}_I(M/\ker(f))={\rm rank}_I(\ker(f))+{\rm rank}_I({\rm im}(f)).$$
\end{proof}

\section{Isomorphic Basis of a Module in the Finite Rank Case}

\begin{definition}
Let $R$ be a commutative ring with identity, and let $M$ be a free $R$-module with ${\rm rank}_I(M)=n$.
The set (or multiset) $\{m_1,\ldots,m_n\}$ for $m_1,\ldots,m_n \in M$ is an {\bf isomorphic basis} of $M$ if and only if the $R$-module homomorphism $f:R^n \to M$ such that $f(e_i)=m_i$ for all $i=1,\ldots,n$ is an isomorphism.
\end{definition}

The standard basis of $R^n$, $\{e_1,\ldots,e_n\}$, is also an isomorphic basis of $R^n$ according to this definition using the identity map.

Let $R$ be a commutative ring with identity, and let $M$ be a module over $R$. Since a module homomorphism with domain $R^n$ is uniquely determined by the outputs of the standard basis, for ${\rm rank}_I(M)=n$, there is a one-to-one correspondence between the module isomorphisms $f:R^n \to M$ and the ordered bases of $M$.

\begin{theorem}
Let $R$ be a commutative ring with identity, and let $M$ be an $R$-module. Then $M$ is free with ${\rm rank}_I(M)=n$ if and only if there exists an isomorphic basis of $M$ containing exactly $n$ elements.
\end{theorem}

\begin{proof}
If ${\rm rank}_I(M)=n$, then there exists a module isomorphism $f:R^n \to M$. Then $\{f(e_1),\ldots,f(e_n)\}$ is an isomorphic basis of $M$ by definition.

If $\{m_1,\ldots,m_n\}$ is an isomorphic basis of $M$, then the module homomorphism $f:R^n \to M$ defined by $f(e_i)=m_i$ for all $i=1,\ldots,n$ is an isomorphism. Since there exists a module isomorphism between $R^n$ and $M$, $M$ is free with ${\rm rank}_I(M)=n$.
\end{proof}

If $M$ is a free module with a finite rank, ${\rm rank}(M)$ is unique, so we get the following corollary.

\begin{corollary}
Let $R$ be a commutative ring with identity, and let $M$ be a free $R$-module with finite rank. Every isomorphic basis of $M$ contains exactly ${\rm rank}_I(M)$ elements.
\end{corollary}

When ${\rm rank}_I(M)=0$, the only isomorphic basis of $M$ is the empty set.

\begin{definition}
Let $R$ be a commutative ring with identity, and let $M$ be an $R$-module.
For the set (or multiset) $\{m_1,\ldots,m_n\}$ with $m_1,\ldots,m_n \in M$, an {\bf $R$-linear combination} of $\{m_1,\ldots,m_n\}$ is an output of the module homomorphism $f:R^n \to M$ such that $f(e_i)=m_i$ for all $i=1,\ldots,n$. 
\end{definition}

This means that an $R$-linear combination is $f(x_1,\ldots,x_n)=x_1m_1+\cdots+x_nm_n$ for some $x_1,\ldots,x_n \in R$, which is the same as the usual definition.

\begin{theorem}
Let $R$ be a commutative ring with identity, let $M$ be an $R$-module, and let $\mathcal{B}=\{m_1,\ldots,m_n\}$ for some $m_1,\ldots,m_n \in M$. Then $\mathcal{B}$ is an isomorphic basis of $M$ if and only if every element in $M$ is a unique linear combination of $\mathcal{B}$.
\end{theorem}

\begin{proof}
Let $f:R^n \to M$ be the module homomorphism such that $f(e_i)=m_i$ for all $i=1,\ldots,n$. 

If $\mathcal{B}$ is an isomorpic basis, then $f$ is an isomorphism and every element in $M$ is a unique output for some input $(x_1,\ldots,x_n) \in R^n$, which means every element in $M$ is a unique $R$-linear combination $x_1m_1+\cdots+x_nm_n$ for that $(x_1,\ldots,x_n) \in R^n$. 

If every element in $M$ is a unique $R$-linear combination of $\mathcal{B}$, then $f$ is surjective because every element in $M$ is an output of $f$, and $f$ is injective because $f(0,\ldots,0)=0v_1+\cdots+0v_n={\bf 0}_M$ is unique, which implies $\ker(f)=\{{\bf 0}\}$. Thus $f$ is an isomorphism, and $\mathcal{B}$ is an isomorphic basis.
\end{proof}

\section{Surjective and Injective Sets in the Finite Rank Case}

\begin{definition}
Let $R$ be a commutative ring with identity, and let $M$ be an $R$-module.
For the set (or multiset) $S=\{m_1,\ldots,m_n\}$ for $m_1,\ldots,m_n \in M$, the {\bf $R$-span} of $S$ is the image of the module homomorphism $f:R^n \to M$ such that $f(e_i)=m_i$ for all $i=1,\ldots,n$, or in other words, ${\rm span}_{R}(S)={\rm im}(f)$. 
\end{definition}

This means that ${\rm span}_{R}(m_1,\ldots,m_n)=\{x_1m_1+\cdots+x_nm_n|x_1,\ldots,x_n \in R\}$, which is the same as the usual definition.

\begin{definition}
Let $R$ be a commutative ring with identity, and let $M$ be an $R$-module.
The set (or multiset) $\{m_1,\ldots,m_n\}$ for $m_1,\ldots,m_n \in M$ is a {\bf surjective set} in $M$ if and only if the module homomorphism $f:R^n \to M$ such that $f(e_i)=m_i$ for all $i=1,\ldots,n$ is surjective.
\end{definition}

Therefore, $S=\{m_1,\ldots,m_n\}$ is a surjective set in $M$ if and only if ${\rm span}_{R}(S)=M$. The module homomorphism $f:R^n \to M$ defined by $f(e_i)=m_i$ for all $i=1,\ldots,n$ is also surjective if and only if $S$ is a generating set of $M$. Therefore the definitions of surjective set in $M$ and generating set of $M$ are equivalent (for finitely-generated $M$ as currently defined). 

\begin{theorem}
Let $R$ be a commutative ring with identity, let $M$ be an $R$-module, and let $S=\{m_1,\ldots,m_n\}$ be a set or multiset of elements in $M$. Then $S$ is a generating set of $M$ if and only if $S$ is a surjective set in $M$.
\end{theorem}

We can prove the results involving generating sets for surjective sets using the definition of surjective sets.

\begin{theorem} \label{surjset>=rank(M)}
Let $R$ be a commutative ring with identity, let $M$ be a free $R$-module with finite rank, and let $\{m_1,\ldots,m_k\}$ be a surjective set in $M$. Then $k \geq {\rm rank}(M)$.
\end{theorem}

\begin{proof}
Since $\{m_1,\ldots,m_k\}$ is a surjective set in $M$, there exists a surjective $R$-module homomorphism $g:R^k \to M$ such that $g(e_i)=m_i$ for all $i=1,\ldots,k$.

Since ${\rm rank}_I(R^k)=k$, $k \geq {\rm rank}(M)$ by Theorem \ref{thm:surj:rank(M)>=rank(N)}.
\end{proof}

\begin{definition}
Let $R$ be a commutative ring with identity, and let $M$ be an $R$-module.
The set (or multiset) $\{m_1,\ldots,m_n\}$ for $m_1,\ldots,m_n \in M$ is an {\bf injective set} in $M$ if and only if the module homomorphism $f:R^n \to M$ such that $f(e_i)=m_i$ for all $i=1,\ldots,n$ is injective.
\end{definition}

If $n=0$ and $\{m_1,\ldots,m_n\}=\emptyset$, then $R^0=\{{\bf 0}\}$ and the module homomorphism $f:\{{\bf 0}\} \to M$ is defined by $f({\bf 0})={\bf 0}_M$, so $f$ is injective. This implies the empty set is an injective set in $M$.

\begin{theorem}
Let $R$ be a commutative ring with identity, let $M$ be an $R$-module, and let $S=\{m_1,\ldots,m_n\}$ be a set or multiset of elements in $M$. Then $S$ is a linearly independent set in $M$ if and only if $S$ is an injective set in $M$.
\end{theorem}

\begin{proof}
Let $f:R^n \to M$ be the module homomorphism defined by $f(e_i)=m_i$ for all $i=1,\ldots,n$. Then $f(x_1,\ldots,x_n)=x_1m_1+\cdots+x_nm_n$ for all $x_1,\ldots,x_n \in R$ and $(x_1,\ldots,x_n) \in \ker(f)$ if and only if $x_1m_1+\cdots+x_nm_n={\bf 0}_M$. 

If $S$ is a linearly independent set in $M$, then $x_1m_1+\cdots+x_nm_n={\bf 0}_M$ implies $(x_1,\ldots,x_n)=(0,\ldots,0)$, $\ker(f)=\{(0,\ldots,0)\}$, and $f$ is injective, so $S$ is an injective set in $M$. If $S$ is an injective set in $M$, then $f$ is injective and $\ker(f)=\{(0,\ldots,0)\}$, and $x_1m_1+\cdots+x_mv_n={\bf 0}_M$ implies $(x_1,\ldots,x_n)=(0,\ldots,0)$.
\end{proof}

\section{Isomorphic Dimension with the Infinite Case}

Let $F$ be a field, and let $S$ be a set. Let $\{e_i|i \in S\}$ be the standard basis of $F$ indexed by $S$, and define $F_0^S$ to be the vector space of all linear combinations of finite subsets of the standard basis of $F$ indexed by $S$. We can also think of $F_0^S$ as the vector space of all functions with finite support from $S$ to $F$, and here $e_i$ is the function from $S$ to $F$ such that $e_i(s)=0$ if $s \neq i$ and $e_i(s)=1$ if $s=i$. If $S$ is finite with $|S|=n$, then we can apply our results for $F^n$ to $F_0^S$ because these are isomorphic by sending the standard basis of $F^n$ to the standard basis of $F$ indexed by $S$.

\begin{theorem}
Let $V$ be a vector space over field $F$,and let $S$ be a set. Then for any set (or multiset) $\{v_i|i \in S\}$ with $v_i \in V$ for all $i \in S$, $f:F_0^S \to V$ such that
$$f(x_{i_1}e_{i_1}+\cdots+x_{i_n}e_{i_n})=x_{i_1}v_{i_1}+\cdots+x_{i_n}v_{i_n}$$
for all nonnegative integers $n$ and all $x_{i_1},\ldots,x_{i_n} \in F$, is the unique linear transformation such that
$f(e_i)=v_i$ for all $i \in S$.
\end{theorem}

When $S=\emptyset$ and $|S|=0$, $F_0^S=\{{\bf 0}\}$, and the only linear transformation $f:\{{\bf 0}\} \to V$ is defined by $f({\bf 0})={\bf 0}_V$.

\begin{lemma} \label{Infinite:lemma:Sfinite}
Let $F$ be a field, let $S$ be a set, and let $n$ be a finite nonnegative integer. $|S|=n$ if and only if $\dim_I(F_0^S)=n$.
\end{lemma}

\begin{proof}
If $S$ is finite with $|S|=n$, then $S=\{s_1,\ldots,s_n\}$ for some distinct $s_1,\ldots,s_n \in S$. The linear transformation $f:F^n \to F_0^S$ such that $f(e_i)=e_{s_i}$ for all $i=1,\ldots,n$ is bijective and thus a linear isomorphism. By definition, $\dim_I(F_0^S)=n$.

If $S$ is infinite, then $\aleph_0 \leq |S|$, and there exists a subset of $S$ that can be written as $\{s_i|i \in \mathbb{N}\}$ of distinct elements from $S$. For every nonnegative integer $n$, define $f_n:F^n \to F_0^S$ as the linear transformation such that $f(e_i)=e_{s_i}$ for all $i=1,\ldots,n$. This creates a sequence of injective linear transformations.
$$f_{0}:F^0 \to F_0^S, \, f_{1}:F^1 \to F_0^S, \ldots, \, f_{n}:F^n \to F_0^S, \ldots$$
such that $f_{k+1} \circ p_k^{k+1} = f_k$ and ${\rm im}(f_{k}) \subsetneq {\rm im}(f_{k+1})$ for any $k \geq 0$.
By Theorem \ref{injseqthm}, $F_0^S$ is infinite-dimensional.
\end{proof}

\begin{lemma} \label{Infinite:injlemma}
Let $F$ be a field, and let $S$ and $T$ be sets. If $f:F_0^S \to F_0^T$ is an injective linear transformation, then $|S| \leq |T|$.
\end{lemma}

\begin{proof}
If $T$ is finite, then by Lemma \ref{Infinite:lemma:Sfinite} and by Theorem \ref{thm:inj:dim(V)<=dim(W)}, $S$ is finite and $|S|=\dim_I(F_0^S) \leq \dim_I(F_0^T)=|T|$.

Consider $T$ to be infinite. For any $i \in S$, $f(e_i)$ is a linear combination of a finite subset of $\{e_j|j \in T\}$. Let $A_i$ be the finite subset of $T$ of the indexes of the basis elements in the linear combination of $f(e_i)$. 

For any finite subset $T'$ of $T$, let $S'=\{i \in S|A_i = T'\}$. Either $S'=\emptyset$ and $|S'|=0$, or if $|S'|>0$, then the linear transformation $f':F_0^{S'} \to F_0^{T'}$ such that $f'(e_i)=f(e_i)$ for all $i \in S'$ is an injective function and $|S'|$ is finite with $|S'| \leq |T'|$.

Define $g:S \to \mathcal{P}_0(T) \times \mathbb{N}$ such that $g(i)=A_i \times j$ for every $i \in S$, and since for every $A \in \mathcal{P}_0(T)$, there will only be a finite number of $i$ with $A_i=A$, we can order these using the positive integers $j$, $1 \leq j \leq |A|$ for those $i$ with $A_i=A$. The function $g$ is injective, and $|S| \leq |\mathcal{P}_0(T) \times \mathbb{N}| = |T|$ because $T$ is infinite.
\end{proof}

\begin{lemma} \label{Infinite:surjlemma}
Let $F$ be a field, and let $S$ and $T$ be sets. If $f:F_0^S \to F_0^T$ is a surjective linear transformation, then $|S| \geq |T|$.
\end{lemma}

\begin{proof}
If $S$ is finite, then by Lemma \ref{Infinite:lemma:Sfinite} and Theorem \ref{thm:surj:dim(V)>=dim(W)}, $|S|=\dim_I(F_0^S) \geq \dim_I(F_0^T)=|T|$.

Consider $S$ to be infinite. For any $i \in S$, $f(e_i)$ is a linear combination of a finite subset of $\{e_j|j \in T\}$. Let $A_i$ be the finite subset of $T$ of the indexes of the basis elements in the linear combination of $f(e_i)$. Since $f$ is surjective, $\cup_{i \in S} A_i = T$. Therefefore $|T| = |\cup_{i \in S} A_i| \leq |S \times \mathbb{N}| = |S| \cdot \aleph_0 = |S|$ because $S$ is infinite.
\end{proof}

\begin{proposition} \label{|S|vs|T|}
Let $V$ be a vector space over field $F$, and let $S$ and $T$ be sets. Let $f_S:F_0^S \to V$ be a linear transformation and let $f_T:F_0^T \to V$ be a linear isomorphism. If $f_S$ is injective, then $|S| \leq |T|$, and if $f_S$ is surjective, then $|S| \geq |T|$.
\end{proposition}

\begin{proof}
If $f_S$ is injective, then $f_T^{-1} \circ f_S:F_0^S \to F_0^T$ is an injective linear transformation, and by Lemma \ref{Infinite:injlemma}, $|S| \leq |T|$.

If $f_S$ is surjective, then $f_T^{-1} \circ f_S:F_0^S \to F_0^T$ is a surjective linear transformation, and by Lemma \ref{Infinite:surjlemma}, $|S| \geq |T|$.
\end{proof}

\begin{corollary} \label{|S|=|T|}
Let $V$ be a vector space over field $F$, and let $S$ and $T$ be sets. If $f_S:F_0^S \to V$ is a linear isomorphism and $f_T:F_0^T \to V$ is a linear isomorphism, then $|S|=|T|$.
\end{corollary}

\begin{proof}
By Proposition \ref{|S|vs|T|}, $|S| \leq |T|$ and $|S| \geq |T|$, so $|S|=|T|$ by the Schr{\" o}der-Bernstein Theorem.
\end{proof}

\begin{theorem}
Let $V$ be a vector space over field $F$, and let $S$ and $T$ be sets. If $f_S:F_0^S \to V$ is a linear isomorphism and $|S|=|T|$, then there exists a linear isomorphism $f_T:F_0^T \to V$.
\end{theorem}

\begin{proof}
If $|S|=|T|$, then there exists an isomorphism $g:T \to S$. Define the linear transformation $f_T:F_0^T \to V$ by $f_T(e_i)=f_S(e_{g(i)})$ for all $i \in T$. Since $g$ is an isomorphism and $f_S$ is a linear isomorphism, then $f_T$ is a linear isomorphism.
\end{proof}

If there exists a linear isomorphism between $F_0^S$ and $V$, then $|S|$ is unique, and for any set $T$ with the cardinality $|T|=|S|$, there exists a linear isomorphism between $F_0^T$ and $V$. We can generalize our previous definition of isomorphic dimension to include dimensions of other cardinalities. This definition coincides with our previous definition because $F^n=F_0^S$ for $S=\{1,\ldots,n\}$, where $\{e_i|i \in S\}=\{e_1,\ldots,e_n\}$ is the standard basis of $F^n$.

\begin{definition}
Let $V$ be a vector space over field $F$, and let $S$ be a set. The {\bf isomorphic dimension} of $V$ is $|S|$, labeled $\dim_{I}(V)=|S|$, if and only if there exists a linear isomorphism $f:F_0^S \to V$.
\end{definition}

For every pair of sets $S \subset T$, define $p_S^T:F_0^S \to F_0^T$ to be the linear transformation such that $p_S^T(e_i)=e_i$ for all if $i \in S$, and define $p_T^S:F_0^T \to F_0^S$ to be the linear transformation such that $p_T^S(e_i)=e_i$ for all $i \in S$ and $p_T^S(e_i)={\bf 0}$ for all $i \in T-S$.

\begin{lemma} \label{infinite:addinjective}
Let $V$ be a vector space over field $F$, and let $S$ be a set. Let $f:F_0^S \to V$ be a linear transformation with $f(e_i)=v_i$ for all $i \in S$.
If $f$ is injective but not surjective, then for any $v \in V$ with $v \notin {\rm im}(f)$ and any $j \notin S$, there exists an injective linear transformation $f':F_0^{S'} \to V$ for $S'=S \cup \{j\}$ such that $f'(e_i)=v_i$ for all $i \in S$ and $f'(e_j)=v$ with $f' \circ p_S^{S'} = f$ and ${\rm im}(f) \subsetneq {\rm im}(f')$.
\end{lemma}

\begin{proof}
If $f:F_0^S \to V$ is not surjective, there exists $v \in V$ such that $v \notin {\rm im}(f)$. Then for any $a \in F$, $av \in {\rm im}(f)$ if and only if $a=0$ because ${\rm im}(f)$ is a subspace.

For $S'=S \cup \{j\}$, define $f':F_0^{S'} \to V$ to be the linear transformation such that $f'(e_i)=v_i$ for all $i \in S$ and $f'(e_j)=v$. It follows that $f' \circ p_S^{S'} = f$ and ${\rm im}(f) \subset {\rm im}(f')$ with ${\rm im}(f) \neq {\rm im}(f')$ because $v \notin {\rm im}(f)$ but $v \in {\rm im}(f')$.

For any $x_1e_{i_1}+\cdots+x_ne_{i_n}+xe_{j} \in \ker(f')$, $$f(x_1e_{i_1}+\cdots+x_ne_{i_n})+xv={\bf 0},$$ which implies $f(x_1e_{i_1}+\cdots+x_ne_{i_n})=-xv$. Therefore, $f(x_1e_{i_1}+\cdots+x_ne_{i_n})={\bf 0}_V$ and $x=0$. Since $f$ is injective, $x_1e_{i_1}+\cdots+x_ne_{i_n}={\bf 0}$, and $x_1e_{i_1}+\cdots+x_ne_{i_n}+xe_{j}={\bf 0}+0e_{j}={\bf 0}$. This proves $f'$ is injective.
\end{proof}

\begin{lemma} \label{lemma:SCBCT}
Let $V$ be a vector space over field $F$, and let $S$ and $T$ be sets such that $S \subset T$. If $f_T:F_0^T \to V$ is a surjective linear transformation such that the linear transformation $f_S:F_0^S \to V$ defined by $f_S(e_i)=f_T(e_i)$ for all $i \in S$ is injective, then there exists a set $B$, such that $S \subset B \subset T$ and the linear transformation $f_B:F_0^B \to V$ defined by $f_B(e_i)=f_T(e_i)$ for all $i \in B$ is an isomorphism.
\end{lemma}

\begin{proof}
Let $U$ be the set of all subsets $A$ such that $S \subset A \subset T$ and the linear transformation $f_A:F_0^A \to V$ defined by $f_A(e_i)=f_T(e_i)$ for all $i \in A$ is injective. Since $S \in U$, $U$ is nonempty.

Let $\{A_j|j \in K\}$ be a chain in $U$. Then $C=\cup_{j \in K} A_j$ is also in $U$ because $S \subset C \subset T$ and the linear transformation $f_C:F_0^C \to V$ defined by $f_C(e_i)=f_T(e_i)$ for all $i \in C$ is injective. Therefore $C$ is an upperbound for chain $\{A_j|j \in K\}$. Since every chain in $U$ has an upper bound, $U$ has a maximal element $B$ by Zorn's Lemma.

Since $B$ is in $U$, $S \subset B \subset T$ and the linear transformation $f_B:F_0^B \to V$ defined by $f_B(e_i)=f_T(e_i)$ for all $i \in B$ is injective. If $f_T(e_i) \in {\rm im}(f_B)$ for all $i \in T$, then $f_B$ is surjective because $f_T$ is surjective. If $f_T(e_j) \notin {\rm im}(f_B)$ for some $j \in T-B$, then for $B'=B \cup \{j\}$, $S \subset B' \subset T$ and the linear transformation $f_{B'}:F_0^{B'} \to V$ defined by $f_{B'}(e_i)=f_T(e_i)$ for all $i \in B'$ is injective by Lemma \ref{infinite:addinjective}. This contradicts the maximality of $B$, so $f_T(e_i) \in {\rm im}(f_B)$ for all $i \in T$ and $f_B$ is surjective.

Therefore $f_B$ is a linear isomorphism.
\end{proof}

\begin{theorem} \label{thm:Bexists}
Let $V$ be a vector space over field $F$. There exists a set $B$ such that $\dim_I(V)=|B|$.
\end{theorem}

\begin{proof}
Let $S = \emptyset$ and $T = V$. Then $S \subset T$. Define the linear transformation $f_T:F_0^T \to V$ by $f_T(e_i)=i$ for all $i \in T$. Then $f_T$ is surjective. The linear transformation $f_S:F_0^S \to V$ defined by $f_S(e_i)=f_T(e_i)$ for all $i \in S$ is the same as the linear transformation $f_S:F_0^S \to V$ defined by $f_S({\bf 0})={\bf 0}_V$, which is injective. By Lemma \ref{lemma:SCBCT}, there exists a set $B$ such that $S \subset B \subset T$ and a linear isomorphism $f_B:F_0^B \to V$. By definition, $\dim_I(V)=|B|$.
\end{proof}

The isomorphic definition of dimension allows for a short proof of the following.

\begin{theorem} \label{Infinite:isoequaldim}
Let $V$ and $W$ be vector spaces over field $F$. There exists linear isomorphism $f:V \to W$ if and only if $\dim_{I}(V)=\dim_{I}(W)$.
\end{theorem}

\begin{proof}
Let $\dim_{I}(V)=|S|$. There exists linear isomorphism $f_V:F_0^S \to V$. 

If $f:V \to W$ is a linear isomorphism, then $f \circ f_V:F_0^S \to W$ is a linear isomorphism. Therefore $\dim_I(W)=|S|$.

If $\dim_I(W)=|S|$, then there exists linear isomorphism $f_W:F_0^S \to W$, and $f_W \circ f_V^{-1}:V \to W$ is a linear isomorphism.
\end{proof}

We also say $V$ and $W$ are isomorphic if and only if there exists a linear isomorphism between them, so $V$ and $W$ are isomorphic if and only if $\dim_I(V)=\dim_I(W)$.

\begin{theorem}
Let $V$ be a vector space over field $F$, and let $U$ be a subspace of $V$. Then $\dim_I(U) \leq \dim_I(V)$.
\end{theorem}

\begin{proof}
By Theorem \ref{thm:Bexists}, there exist a set $S$ such that $\dim_I(U)=|S|$ and set $T$ such that $\dim_I(V)=|T|$. By definition, there exist linear isomorphisms $f_S:F_0^S \to U$ and $f_T:F_0^T \to V$. Let $g:U \to V$ be the injective linear transformation defined by $g(u)=u$ for all $u \in U$. Then $f_T^{-1} \circ g \circ f_S:F_0^S \to F_0^T$ is injective. By Lemma \ref{Infinite:injlemma}, $\dim_I(U)=|S| \leq |T|=\dim_I(V)$.
\end{proof}

\begin{theorem} \label{Infinite:thm:inj:dim(V)<=dim(W)}
Let $V$ and $W$ be vector spaces over field $F$. There exists an injective linear transformation $f:V \to W$ if and only if $\dim_I(V) \leq \dim_I(W)$.
\end{theorem}

\begin{proof}
By Theorem \ref{thm:Bexists}, there exist a set $S$ such that $\dim_I(V)=|S|$ and a set $T$ such that $\dim_I(W)=|T|$. By definition, there exist linear isomorphisms $f_V:F_0^S \to V$ and $f_W:F_0^T \to W$.

If $f:V \to W$ is an injective linear transformation, then $f_W^{-1} \circ f \circ f_V:F_0^S \to F_0^T$ is an injective linear transformation, and by Lemma \ref{Infinite:injlemma}, $\dim_I(V) = |S| \leq |T| = \dim_I(W)$.

If $|S| = \dim_I(V) \leq \dim_I(W) = |T|$, then there exists an injective function $g:S \to T$. Define the linear transformation $g':F_0^S \to F_0^T$ such that $g'(e_i)=e_{g(i)}$ for all $i \in S$. Since $g$ is injective, $g'$ is an injective linear transformation. The linear transformation $f_W \circ g' \circ f_V^{-1}:V \to W$ is injective.
\end{proof}

\begin{theorem} \label{Infinite:thm:surj:dim(V)>=dim(W)}
Let $V$ and $W$ be vector spaces over field $F$. There exists a surjective linear transformation $f:V \to W$ if and only if $\dim_I(V) \geq \dim_I(W)$.
\end{theorem}

\begin{proof}
By Theorem \ref{thm:Bexists}, there exist a set $S$ such that $\dim_I(V)=|S|$ and a set $T$ such that $\dim_I(W)=|T|$. By definition, there exist linear isomorphisms $f_V:F_0^S \to V$ and $f_W:F_0^T \to W$.

If $f:V \to W$ is a surjective linear transformation, then $f_W^{-1} \circ f \circ f_V:F_0^S \to F_0^T$ is an surjective linear transformation, and by Lemma \ref{Infinite:surjlemma}, $\dim_I(V) = |S| \geq |T| = \dim_I(W)$.

If $|S| = \dim_I(V) \geq \dim_I(W) = |T|$, then there exists a surjective function $g:S \to T$. Define the linear transformation $g':F_0^S \to F_0^T$ such that $g'(e_i)=e_{g(i)}$ for all $i \in S$. Since $g$ is surjective, $g'$ is a surjective linear transformation. The linear transformation $f_W \circ g' \circ f_V^{-1}:V \to W$ is surjective.
\end{proof}

\begin{theorem} \label{Infinite:dim(V/U)}
Let $V$ be a vector space over field $F$, and let $U$ be a subspace of $V$. 

$$\dim_I(U)+\dim_I(V/U)=\dim_I(V)$$
\end{theorem}

\begin{proof}
By Theorem \ref{thm:Bexists}, there exists a set $S$ such that $\dim_I(U)=|S|$, and there exists a linear isomorphism $f:F_0^S \to U$. Let ${\rm inc}:U \to V$ be the injective linear transformation such that ${\rm inc}(u)=u$ for any $u \in U$, and let $f_S={\rm inc} \circ f:F_0^S \to V$. As a composition of injective linear transformations, $f_S$ is injective.

Let $T=S \sqcup (V-\{f(e_i)|i \in S\})$. Then the linear transformation $f_T:F_0^T \to V$ defined as $f_T(e_i) = f_S(e_i)$ for all $i \in S$ and $f_T(e_i) = i$ for all $i \in T-S$ is surjective. By Lemma \ref{lemma:SCBCT}, there exist a set $B$ such that $S \subset B \subset T$ and a linear isomorphism $g:F_0^B \to V$ such that $g(e_i)=f_S(e_i)$ for all $i \in S$ and thus ${\rm im}(g|_{F_0^S})=U$. Since $g$ is a linear isomorphism, $\dim_I(V)=|B|$.

Let $q:V \to V/U$ be the surjective linear transformation such that $q(v)=v+U$ for all $v \in V$. Then $h=q \circ g|_{F_0^{B-S}}:F_0^{B-S} \to V/U$ is a linear transformation. Since $g$ is a injective, ${\rm im}(g|_{F_0^S}) \cap {\rm im}(g|_{F_0^{B-S}})=\{{\bf 0}\}$. Thus $U \cap {\rm im}(g|_{F_0^{B-S}})=\{{\bf 0}\}$ and $h$ is injective. Since $g$ is surjective, for any $v+U$ in $V/U$, $v$ is a linear combination of the outputs from $g$ of the standard basis of $F_0^B$. This linear combination can be split into a sum of a linear combination with standard basis vectors from $F_0^S$ and and a linear combination from $F_0^{B-S}$ because $B=S \sqcup (B-S)$. Any linear combination of outputs from $g$ of standard basis vectors from $F_0^S$ is in $U$. Thus $v+U=v'+U$ for some $v'$ equal to a linear combination of outputs from $g$ of standard basis vectors from $F_0^{B-S}$, and $h$ is surjective.

As a bijective linear transformation, $h$ is a linear isomorphism and $\dim_I(V/U)=|B-S|$.
Therefore $|S|+|B-S|=|B|$ or $\dim_I(U)+\dim_I(V/U)=\dim_I(V)$.
\end{proof}

The following theorem is the Rank Nullity Theorem using the definition of isomorphic dimension.

\begin{theorem}
Let $V$ and $W$ be vector spaces over field $F$, and let $f:V \to W$ be a linear transformation.

$$\dim_I(\ker(f))+\dim_I({\rm im}(f))=\dim_I(V)$$
\end{theorem}

\begin{proof}
By the first isomorphism theorem, $\bar{f}:V/\ker(f) \to {\rm im}(f)$ is a linear isomorphism, and by Theorem \ref{Infinite:isoequaldim}, $\dim_I(V/\ker(f))=\dim_I({\rm im}(f))$. By Theorem \ref{Infinite:dim(V/U)}, $\dim_I(\ker(f))+\dim_I({\rm im}(f))=\dim_I(V)$.
\end{proof}

\section{Isomorphic Basis of a Vector Space with the Infinite Dimension Case}

In a similar fashion to our latest definition of dimension, we can also construct a new definition for basis that is equivalent to the algebraic definition of a basis.

\begin{definition}
Let $V$ be a vector space over field $F$, and let $S$ be a set. The set (or multiset) $\{v_i|i \in S\}$ for $v_i \in V$ for all $i \in S$ is an {\bf isomorphic basis} of $V$ if and only if the linear transformation $f:F_0^S \to V$ such that $f(e_i)=v_i$ for all $i \in S$ is an isomorphism.
\end{definition}

The standard basis of $F_0^S$, $\{e_i|i \in S\}$, is also an isomorphic basis of $F_0^S$ according to this definition using the identity map $I:F_0^S \to F_0^S$, which is the linear transformation with $I(e_i)=e_i$ for all $i \in S$ and is an isomorphism.

Let $V$ be a vector space over field $F$. Since a linear transformation with domain $F_0^S$ is uniquely determined by the outputs of the standard basis, for $\dim_I(V)=|S|$, there is a one-to-one correspondence between the linear isomorphisms $f:F_0^S \to V$ and the ordered bases of $V$ (based on the ordering of $S$).

\begin{theorem}
Let $V$ be a vector space over field $F$, and let $S$ be a set. Then $\dim_I(V)=|S|$ if and only if there exists an isomorphic basis of $V$ containing $|S|$ vectors.
\end{theorem}

\begin{proof}
If $\dim_I(V)=|S|$, then there exists a linear isomorphism $f:F_0^S \to V$. Then $\{f(e_i)|i \in S\}$ is an isomorphic basis of $V$ by definition.

If $\{v_i|i \in S\}$ is an isomorphic basis of $V$, then the linear transformation $f:F_0^S \to V$ defined by $f(e_i)=v_i$ for all $i \in S$ is an isomorphism. Since there exists an isomorphism between $F_0^S$ and $V$, $\dim_I(V)=|S|$.
\end{proof}

For vector space $V$, $\dim_I(V)$ is unique, so we get the following corollary.

\begin{corollary}
Let $V$ be a vector space over field $F$. Every isomorphic basis of $V$ has cardinality equal to $\dim_I(V)$.
\end{corollary}

\begin{definition}
Let $V$ be a vector space over field $F$, and let $S$ be a set. For the set (or multiset) $\{v_i|i \in S\}$ with $v_i \in V$ for all $i \in S$, a {\bf linear combination} of $\{v_i|i \in S\}$ is an output of the linear transformation $f:F_0^{S'} \to V$ such that $f(e_i)=v_i$ for all $i \in S'$ for a finite subset $S' \subset S$. 
\end{definition}

This means that a linear combination is $f(x_1e_{i_1}+\cdots+x_ne_{i_n})=x_1v_{i_1}+\cdots+x_nv_{i_n}$ for some nonnegative integer $n$ and some $x_1,\ldots,x_n \in F$, which is the same as the usual definition.

\begin{theorem}
Let $V$ be a vector space over field $F$, let $S$ be a set, and let $\mathcal{B}=\{v_i|i \in S\}$ such that $v_i \in V$ for all $i \in S$. Then $\mathcal{B}$ is an isomorphic basis of $V$ if and only if every vector in $V$ is a unique linear combination of $\mathcal{B}$.
\end{theorem}

\begin{proof}
Let $f:F_0^S \to V$ be the linear transformation such that $f(e_i)=v_i$ for all $i \in S$.

If $\mathcal{B}$ is an isomorphic basis, then $f$ is an isomorphism and every vector in $V$ is a unique output for some input $x_1e_{i_1}+\cdots+x_ne_{i_n} \in F_0^S$, which means it is a unique linear combination $x_1v_{i_1}+\cdots+x_nv_{i_n}$.

If every vector in $V$ is a unique linear combination of $\mathcal{B}$, then $f$ is surjective because every vector in $V$ is an output of $f$, and $f$ is injective because $f({\bf 0})={\bf 0}_V$ is unique, which implies $\ker(f)=\{{\bf 0}\}$. Thus $f$ is an isomorphism, and $\mathcal{B}$ is an isomorphic basis.
\end{proof}

\section{Surjective and Injective Sets with the Infinite Dimension Case}

\begin{definition}
Let $V$ be a vector space over field $F$, and let $S$ be a set. For the set (or multiset) $\{v_i|i \in S\}$ with $v_i \in V$ for all $i \in S$, the {\bf span} of $\{v_i|i \in S\}$ is the image of the linear transformation $f:F_0^S \to V$ such that $f(e_i)=v_i$ for all $i \in S$, or in other words, ${\rm span}(\{v_i|i \in S\})={\rm im}(f)$. 
\end{definition}

This means that ${\rm span}(\{v_i|i \in S\})$ is the set of all linear combinations of finite subsets of $\{v_i|i \in S\}$, which is the same as the usual definition.

\begin{definition}
Let $V$ be a vector space over field $F$, and let $S$ be a set. The set (or multiset) $\{v_i|i \in S\}$ with $v_i \in V$ for all $i \in S$ is a {\bf surjective set} in $V$ if and only if the linear transformation $f:F_0^S \to V$ such that $f(e_i)=v_i$ for all $i \in S$ is surjective.
\end{definition}

Therefore, $\{v_i|i \in S\}$ is a surjective set in $V$ if and only if ${\rm span}(\{v_i|i \in S\})=V$. The linear transformation $f:F_0^S \to V$ defined by $f(e_i)=v_i$ for all $i \in S$ is also surjective if and only if $\{v_i|i \in S\}$ is a spanning set of $V$. Therefore the definitions of surjective set in $V$ and spanning set of $V$ are equivalent. 

\begin{theorem}
Let $V$ be a vector space over field $F$, let $S$ be a set, and let $\{v_i|i \in S\}$ be a set or multiset with $v_i \in V$ for all $i \in S$. Then $\{v_i|i \in S\}$ is a spanning set of $V$ if and only if $\{v_i|i \in S\}$ is a surjective set in $V$.
\end{theorem}

We can prove the results involving spanning sets for surjective sets using the definition of surjective sets.

\begin{theorem} \label{Infinite:surjset>=dim(V)}
Let $V$ be a vector space over field $F$, let $S$ be a set, and let $\{v_i|i \in S\}$ be a surjective set in $V$. Then $|S| \geq \dim_I(V)$.
\end{theorem}

\begin{proof}
By definition, the linear transformation $f_S:F_0^S \to V$ such that $f_S(e_i)=v_i$ for all $i \in S$ is surjective. By Theorem \ref{thm:Bexists}, there exists a set $B$ such that $\dim_I(V)=|B|$, and by definition, there exists a linear isomorphism $f_B:F_0^B \to V$. The linear transformation $f_B^{-1} \circ f_S:F_0^S \to F_0^B$ is surjective. By Lemma \ref{Infinite:surjlemma}, $|S| \geq |B|=\dim_I(V)$.
\end{proof}

\begin{theorem}
Let $V$ be a vector space over field $F$, and let $S$ be a set. If $\{v_i|i \in S\}$ is a surjective set in $V$, then there exists some subset $\{v_i|i \in S'\}$ for some $S' \subset S$, such that $\{v_i|i \in S'\}$ is an isomorphic basis of $V$.
\end{theorem}

\begin{proof}
By definition, the linear transformation $f_S:F_0^S \to V$ such that $f_S(e_i)=v_i$ for all $i \in S$ is surjective. Let $R=\emptyset$. Then $R \subset S$ and the linear transformation $f_R:F_0^R \to V$ defined by $f_R(e_i)=f_S(e_i)$ for all $i \in R$ is the same as the linear transformation $f_R:F_0^R \to V$ defined by $f_R({\bf 0})={\bf 0}_V$, which is injective. By Lemma \ref{lemma:SCBCT}, there exists a set $S'$ such that $R \subset S' \subset S$ and a linear isomorphism $f_{S'}:F_0^{S'} \to V$. By definition, $\{v_i|i \in S'\}$ is an isomorphic basis of $V$.
\end{proof}

We don't have some of the same results for surjective sets of finite-dimensional vector spaces when including infinite-dimensional vector spaces. For example, let $V=P(F)$, the vector space of all polynomials over field $F$. Then $\dim_I(V)=\aleph_0=|S|$ for $S=\mathbb{Z}_{\geq 0}$ because the linear transformation $f:F_0^S \to V$ defined by $f(e_i)=x^i$ for all $i \in S$ is an isomorphism. The linear transformation $g:F_0^S \to V$ defined by $g(e_0)=1$ and $g(e_i)=x^{i-1}$ for all $i \in \mathbb{Z}_{>0}$ is surjective but not an isomorphism.

\begin{definition}
Let $V$ be a vector space over field $F$, and let $S$ be a set. The set (or multiset) $\{v_i|i \in S\}$ with $v_i \in V$ for all $i \in S$ is an {\bf injective set} in $V$ if and only if the linear transformation $f:F_0^S \to V$ such that $f(e_i)=v_i$ for all $i \in S$ is injective.
\end{definition}

If $S=\emptyset$ and $\{v_i|i \in S\}=\emptyset$, then $F_0^S=\{{\bf 0}\}$ and the linear transformation $f:\{{\bf 0}\} \to V$ is defined by $f({\bf 0})={\bf 0}_V$, so $f$ is injective. This implies the empty set is an injective set in $V$.

\begin{theorem}
Let $V$ be a vector space over a field $F$, let $S$ be a set, and let $\{v_i|i \in S\}$ be a set or multiset with $v_i \in V$ for all $i \in S$. Then $\{v_i|i \in S\}$ is a linearly independent set in $V$ if and only if $\{v_i|i \in S\}$ is an injective set in $V$.
\end{theorem}

\begin{proof}
Let $f:F_0^S \to V$ be the linear transformation defined by $f(e_i)=v_i$ for all $i \in S$. Then $f(x_1e_{i_1}+\cdots+x_ne_{i_n})=x_1v_{i_1}+\cdots+x_nv_{i_n}$ for all $x_1e_{i_1}+\cdots+x_ne_{i_n} \in F_0^S$ and $x_1e_{i_1}+\cdots+x_ne_{i_n} \in \ker(f)$ if and only if $x_1v_{i_1}+\cdots+x_nv_{i_n}={\bf 0}_V$. 

If $\{v_i|i \in S\}$ is a linearly independent set in $V$, then $x_1v_{i_1}+\cdots+x_nv_{i_n}={\bf 0}_V$ implies $x_1e_{i_1}+\cdots+x_ne_{i_n}={\bf 0}$, $\ker(f)=\{{\bf 0}\}$, and $f$ is injective, so $\{v_i|i \in S\}$ is an injective set in $V$. If $\{v_i|i \in S\}$ is an injective set in $V$, then $f$ is injective and $\ker(f)=\{{\bf 0}\}$, and $x_1v_{i_1}+\cdots+x_nv_{i_n}={\bf 0}_V$ implies $x_1e_{i_1}+\cdots+x_ne_{i_n}={\bf 0}$.
\end{proof}

We can prove the results involving linearly independent set for injective sets using the definition of injective sets.

\begin{theorem} \label{Infinite:injset<=dim(V)}
Let $V$ be a vector space over field $F$, and let $S$ be a set. If $\{v_i|i \in S\}$ is an injective set in $V$, then $|S| \leq \dim_I(V)$.
\end{theorem}

\begin{proof}
By definition, the linear transformation $f_S:F_0^S \to V$ such that $f_S(e_i)=v_i$ for all $i \in S$ is injective. By Theorem \ref{thm:Bexists}, there exists a set $B$ such that $\dim_I(V)=|B|$, and by definition, there exists a linear isomorphism $f_B:F_0^B \to V$. The linear transformation $f_B^{-1} \circ f_S:F_0^S \to F_0^B$ is injective. By Lemma \ref{Infinite:injlemma}, $|S| \leq |B|=\dim_I(V)$.
\end{proof}

\begin{corollary}
Let $V$ be a vector space over field $F$, and let $S$ and $T$ be sets. If $\{v_i|i \in S\}$ is an injective linearly independent set in $V$ and $\{w_j|j \in T\}$ is a surjective set in $V$, then $|S| \leq |T|$.
\end{corollary}

\begin{proof}
By Theorem \ref{Infinite:surjset>=dim(V)}, $|T| \geq \dim_I(V)$. By Theorem \ref{Infinite:injset<=dim(V)}, $|S| \leq \dim_I(V)$. Therefore $|S| \leq |T|$.
\end{proof}

\begin{theorem}
Let $V$ be a vector space over field $F$, and let $S$ be a set. If $\{v_i|i \in S\}$ is an injective set in $V$, then there exist a set $S'$ with $S \subset S'$ and a set $\{v_i|i \in S'-S\} \subset V$ such that $\{v_i|i \in S'\}$ is an isomorphic basis of $V$.
\end{theorem}

\begin{proof}
By definition, the linear transformation $f_S:F_0^S \to V$ such that $f_S(e_i)=v_i$ for all $i \in S$ is injective. Let $T=S \sqcup (V-\{v_i|i \in S\})$, so that $S \subset T$. Define the linear transformation $f_T:F_0^S \to V$ by $f_T(e_i)=v_i$ for all $i \in S$ and $f_T(e_{v})=v$ for all $v \in V-\{v_i|i \in S\}$, which is surjective. By Lemma \ref{lemma:SCBCT}, there exists a set $S'$ such that $S \subset S' \subset T$ and $f_{S'}:F_0^{S'} \to V$ is a linear isomorphism. Then $\{f_T(e_i)=v_i|i \in S'\}$ is an isomorphic basis.
\end{proof}

We don't have some of the same results for injective sets of finite-dimensional vector spaces when including infinite-dimensional vector spaces. For example, let $V=P(F)$, the vector space of all polynomials over field $F$. Then $\dim_I(V)=\aleph_0=|S|$ for $S=\mathbb{Z}_{\geq 0}$ because the linear transformation $f:F_0^S \to V$ defined by $f(e_i)=x^i$ for all $i \in S$ is an isomorphism. The linear transformation $g:F_0^S \to V$ defined by $g(e_i)=x^{i+1}$ for all $i \in S$ is injective but not an isomorphism.

\section{Isomorphic Rank with the Infinite Case}

Let $R$ be a commutative ring with identity, and let $S$ be a set. Let $\{e_i|i \in S\}$ be the standard basis of $R$ indexed by $S$, and define $R_0^S$ to be the $R$-module of all $R$-linear combinations of finite subsets of the standard basis of $R$ indexed by $S$. We can also think of $R_0^S$ as the $R$-module of all functions with finite support from $S$ to $R$, and here $e_i$ represents the function from $S$ to $R$ such that $e_i(s)=0$ if $s \neq i$ and $e_i(s)=1$ if $s=i$. If $S$ is finite with $|S|=n$, then we can apply our results for $R^n$ to $R_0^S$ because these are isomorphic using the $R$-module homomorphism that sends the standard basis of $R^n$ to the standard basis of $R$ indexed by $S$.

Let $I$ be an ideal of $R$. Define $I_0^S$ to be the subset of $R_0^S$ containing all $I$-linear combinations of finite subsets of the standard basis of $R$ indexed by $S$. We can also think of $I_0^S$ as the subset of $R_0^S$ of all functions with finite support from $S$ to $R$ in which the outputs are restricted to $I$.

\begin{theorem}
Let $R$ be a commutative ring with identity, let $M$ be an $R$-module,and let $S$ be a set. Then for any set (or multiset) $\{m_i|i \in S\}$ with $m_i \in M$ for all $i \in S$, $f:R_0^S \to M$ such that
$$f(x_{i_1}e_{i_1}+\cdots+x_{i_n}e_{i_n})=x_{i_1}m_{i_1}+\cdots+x_{i_n}m_{i_n}$$
for all nonnegative integers $n$ and all $x_{i_1},\ldots,x_{i_n} \in R$, is the unique $R$-module homomorphism such that
$f(e_i)=m_i$ for all $i \in S$.
\end{theorem}

When $S=\emptyset$ and $|S|=0$, $R_0^S=\{{\bf 0}\}$, and the only $R$-module homomorphism $f:\{{\bf 0}\} \to M$ is defined by $f({\bf 0})={\bf 0}_M$.

\begin{theorem} \label{infinite:thmhatf}
Let $R$ be a commutative ring with identity, let $I$ be an ideal of $R$, let $M$ be an $R$-module, and let $S$ be a set. For any surjective $R$-module homomorphism $f:R_0^S \to M$, $$\hat{f}:R_0^S/I_0^S \to M/f(I_0^S)$$ such that $$\hat{f}(r_{j_1}e_{j_1}+\cdots+r_{j_n}e_{j_n}+I_0^S)=f(r_{j_1}e_{j_1}+\cdots+r_{j_n}e_{j_n})+f(I_0^S)$$ for all nonnegative integers $n$ and all $r_{j_1}e_{j_1}+\cdots+r_{j_n}e_{j_n} \in R_0^S$, is a surjective $R/I$-module homomorphism, and if $f$ is bijective, then $\hat{f}$ is bijective.
\end{theorem}

If $S$ is the empty set, all modules are zero modules, and all homomorphisms are defined by sending ${\bf 0}$ to ${\bf 0}$ for the appropriate modules.

\begin{proof}
$I_0^S$ is an $R$-submodule of $R_0^S$ because $I$ is an ideal of $R$, and $f(I_0^S)$ is an $R$-submodule of $M$ because it is the image of an $R$-module homomorphism of an $R$-submodule.

$R_0^S/I_0^S$ is an $R/I$-module with scalar multiplication defined as $$(a+I) \cdot (r_{j_1}e_{j_1}+\cdots+r_{j_n}e_{j_n}+I_0^S) = a(r_{j_1}e_{j_1}+\cdots+r_{j_n}e_{j_n})+I_0^S.$$ Scalar multiplication is well-defined because $I$ is an ideal or $R$.

$M/f(I_0^S)$ is an $R/I$-module with scalar multiplication defined as $$(a+I) \cdot (m+f(I_0^S)) = am+f(I_0^S).$$ Scalar multiplication is well-defined because $I$ is an ideal of $R$ and $f$ is a surjective $R$-module homomorphism.

$\hat{f}$ is well-defined because $f$ is an $R$-module homomorphism and $I_0^S$ is an $R$-submodule.
Since $f$ is an $R$-module homomorphism, it follows directly that $\hat{f}$ is an $R$-module homomorphism and an $R/I$-module homomorphism.

Since $f$ is surjective, for any $m+f(I_0^S) \in M/f(I_0^S)$ with $m \in M$, then there exists $r_{j_1}e_{j_1}+\cdots+r_{j_n}e_{j_n} \in R_0^S$ such that $f(r_{j_1}e_{j_1}+\cdots+r_{j_n}e_{j_n})=m$. Therefore $\hat{f}(r_{j_1}e_{j_1}+\cdots+r_{j_n}e_{j_n}+I_0^S)=m+f(I_0^S)$, and $\hat{f}$ is surjective.

If $r_{j_1}e_{j_1}+\cdots+r_{j_n}e_{j_n}+I_0^S \in \ker(\hat{f})$, then $f(r_{j_1}e_{j_1}+\cdots+r_{j_n}e_{j_n})+f(I_0^S)={\bf 0}_M+f(I_0^S)$ and $f(r_{j_1}e_{j_1}+\cdots+r_{j_n}e_{j_n}) \in f(I_0^S)$. If $f$ is injective, $r_{j_1}e_{j_1}+\cdots+r_{j_n}e_{j_n} \in I_0^S$ and $\ker(\hat{f})=\{I_0^S\}$. Therefore $\hat{f}$ is injective when $f$ is injective.

If $f$ is bijective, then $\hat{f}$ is bijective by the previous two arguments.
\end{proof}

For any $R$-module homomorphism $f$, the definition of $\hat{f}$ in the previous theorem is an $R$-module homomorphism, but we are particularly interested in the case where $f$ is surjective and $\hat{f}$ is an $R/I$-module homomorphism.

For surjective $R$-module homomorphism $f:R_0^S \to M$, if $f|_{I_0^S}$ is injective and $\hat{f}$ is an isomorphism, then $f$ is an isomorphism. 

If $I=\{{\bf 0}\}$, then $I_0^S=\{{\bf 0}\}$ and $f$ and $\hat{f}$ are equivalent. If $R$ is a field, then the only proper ideal of $R$ is $I=\{{\bf 0}\}$, and $f$ and $\hat{f}$ are equivalent linear transformations.

In the following corollary, in which $R$ is a commutative ring with identity, $I$ is an ideal of $R$, and $M$ is an $R$-module, we use $$IM = \{i_1m_1+\cdots+i_km_k|k \in \mathbb{Z}_{\geq 0}, \, i_1,\ldots,i_k \in I, \, m_1,\ldots,m_k \in M\},$$
the set of all finite $I$-linear combinations of elements of module $M$. $IM$ is a submodule of $M$.

\begin{corollary} \label{infinite:corhatf}
Let $R$ be a commutative ring with identity, let $I$ be an ideal of $R$, let $M$ be an $R$-module, and let $S$ be a set. For any surjective $R$-module homomorphism $f:R_0^S \to M$, $$\hat{f}:(R/I)_0^S \to M/IM$$ such that $$\hat{f}((r_{j_1}+I)e_{j_1}+\cdots+(r_{j_n}+I)e_{j_n})=f(r_{j_1}e_{j_1}+\cdots+r_{j_n}e_{j_n})+IM$$ for all nonnegative integers $n$ and all $r_{j_1}e_{j_1}+\cdots+r_{j_n}e_{j_n} \in R_0^S$ is a surjective $R/I$-module homomorphism, and if $f$ is bijective, then $\hat{f}$ is bijective.
\end{corollary}

If $S$ is the empty set, all modules are zero modules, and all homomorphisms are defined by sending ${\bf 0}$ to ${\bf 0}$ for the appropriate modules.

\begin{proof}
Since $I$ is an ideal of $R$ and $f$ is a surjective $R$-homomorphism, $f(I_0^S)=IM$.

$(R/I)_0^S$ and $R_0^S/I_0^S$ are isomorphic $R$-modules and isomorphic $R/I$-modules using the following module isomorphism: $((r_{j_1}+I)e_{j_1}+\cdots+(r_{j_n}+I)e_{j_n}) \mapsto (r_{j_1}e_{j_1}+\cdots+r_{j_n}e_{j_n})+I_0^S$. Compose this isomorphism with $\hat{f}$ from the previous theorem and replace $f(I_0^S)$ with $IM$.
\end{proof}

\begin{theorem} \label{infinite:iso:|S|=|T|}
Let $R$ be a commutative ring with identity, let $M$ be an $R$-module, and let $S$ and $T$ be sets. If there exist an $R$-module isomorphism $f:R_0^S \to M$ and an $R$-module isomorphism $g:R_0^T \to M$, then $|S|=|T|$.
\end{theorem}

\begin{proof}
There exists a maximal ideal $I$ in $R$, and $R/I$ is a field. By the previous corollary, there exist $R/I$-module isomorphisms $\hat{f}:(R/I)_0^S \to M/IM$ and $\hat{g}:(R/I)_0^T \to M/IM$. Since $R/I$ is a field, $\hat{f}$ and $\hat{g}$ are linear isomorphisms between vector spaces. Two vector spaces are isomorphic if and only if they have equal dimensions, so $\dim(R/I)_0^S=\dim(M/IM)=\dim(R/I)_0^T$ or $|S|=|T|$.
\end{proof}

\begin{theorem}
Let $R$ be a commutative ring with identity, let $M$ be an $R$-module, and let $S$ and $T$ be sets. If $f_S:R_0^S \to M$ is an $R$-module isomorphism and $|S|=|T|$, then there exists an $R$-module isomorphism $f_T:R_0^T \to M$.
\end{theorem}

\begin{proof}
If $|S|=|T|$, then there exists an isomorphism $g:T \to S$. Define the $R$-module homomorphism $f_T:R_0^T \to M$ by $f_T(e_t)=f_S(e_{g(t)})$ for all $t \in T$. Since $g$ is an isomorphism and $f_S$ is an $R$-module isomorphism, then $f_T$ is an $R$-module isomorphism.
\end{proof}

If there exists an $R$-module isomorphism between $R_0^S$ and $M$, then $|S|$ is unique, and for any set $T$ with the cardinality $|T|=|S|$, there exists an $R$-module isomorphism between $R_0^T$ and $M$. We can generalize our previous definition of isomorphic rank in the finite case to include ranks of other cardinalities. This definition coincides with our previous definition because $R^n=R_0^S$ for $S=\{1,\ldots,n\}$, where $\{e_i|i \in S\}=\{e_1,\ldots,e_n\}$ is the standard basis of $R^n$.

\begin{definition}
Let $R$ be a commutative ring with identity, let $M$ be an $R$-module, and let $S$ be a set. $M$ is {\bf free} and the {\bf isomorphic rank} of $M$ is $|S|$, labeled ${\rm rank}_{I}(M)=|S|$, if and only if there exists an $R$-module isomorphism $f:R_0^S \to M$.
\end{definition}

The isomorphic definition of rank allows for a short proof of the following.

\begin{theorem} \label{infinite:isoequalrank}
Let $R$ be a commutative ring with identity, let $M$ be a free $R$-module, and let $N$ be an $R$-module. $M$ and $N$ are isomorphic if and only if $N$ is free and ${\rm rank}_I(M)={\rm rank}_I(N)$.
\end{theorem}

\begin{proof}
Since $M$ is free, there exists some set $S$ such that ${\rm rank}_I(M)=|S|$, and there exists an $R$-module isomorphism $f_M:R_0^S \to M$. 

If $M$ and $N$ are isomorphic, there exists an $R$-module isomorphism $g:M \to N$, and $g \circ f_M:R_0^S \to N$ is an $R$-module isomorphism. Therefore $N$ is free and ${\rm rank}_I(N)=|S|$.

If $N$ is free and ${\rm rank}_I(N)=|S|$, then there exists an $R$-module isomorphism $f_N:R_0^S \to N$, and $f_N \circ f_M^{-1}:M \to N$ is an $R$-module isomorphism.
\end{proof}

\begin{theorem} \label{infinite:thm:surj:rank(M)>=rank(N)}
Let $R$ be a commutative ring with identity and let $M$ and $N$ be free $R$-modules. There exists a surjective $R$-module homomorphism $f:M \to N$ if and only if ${\rm rank}_I(M) \geq {\rm rank}_I(N)$.
\end{theorem}

\begin{proof}
Since $M$ and $N$ are free, there exist sets $S$ and $T$ such that ${\rm rank}_I(M)=|S|$ and ${\rm rank}_I(N)=|T|$, and there exist $R$-module isomorphisms $g_M:R_0^S \to M$ and $g_N:R_0^T \to N$.

If there exists a surjective $R$-module homomorphism $f:M \to N$, then $h=g_N^{-1} \circ f \circ g_M:R_0^S \to R_0^T$ is a surjective $R$-module homomorphism.

Let $I$ be a maximal ideal of $R$, so that $R/I$ is a field. Then by Corollary \ref{infinite:corhatf},
there exists a surjective linear transformation $\hat{h}:(R/I)_0^S \to R_0^T/IR_0^T$. Since $I$ is an ideal of $R$, $IR_0^T=I_0^T$. The composition of $\hat{h}$ with the $R/I$-linear isomorphism from $R_0^T/I_0^T$ to $(R/I)_0^T$ defined by $(r_{j_1}e_{j_1}+\cdots+r_{j_n}e_{j_n})+I_0^T \mapsto ((r_{j_1}+I)e_{j_1}+\cdots+(r_{j_n}+I)e_{j_n})$ produces a surjective linear transformation $(R/I)_0^S \to (R/I)_0^T$ over field $R/I$. Therefore $\dim((R/I)_0^S) \geq \dim((R/I)_0^T)$ or $|S| \geq |T|$.

If $|S| \geq |T|$, then there exists a surjective function $q:S \to T$. Define an $R$-module homomorphism $p:R_0^S \to R_0^T$ such that $p(e_s)=e_{q(s)}$ for all $s \in S$. Since $q$ is surjective, $p$ is surjective, and $f=g_N \circ p \circ g_M^{-1}:M \to N$ is a surjective $R$-module homomorphism.
\end{proof}

\begin{theorem} \label{thm:infinite:rank(M/N)}
Let $R$ be a commutative ring with identity, let $M$ be an $R$-module, and let $N$ be a submodule of $M$. If $N$ and $M/N$ are free, then $M$ is free and 
$${\rm rank}_I(M)={\rm rank}_I(N)+{\rm rank}_I(M/N).$$
\end{theorem}

\begin{proof}
Let $S$ and $T$ be sets such that ${\rm rank}_I(N)=|S|$ and ${\rm rank}_I(M/N)=|T|$. There exist $R$-module isomorphisms $f_S:R_0^S \to N$ and $f_T:R_0^T \to M/N$. Define $g:M/N \to M$ to be an injective function that sends each element of $M/N$ to a particular coset representative. Define $f:R_0^{S \sqcup T} \to M$ as the $R$-module homomorphism such that $f(e_i)=f_S(e_i)$ if $i \in S$ and $f(e_i)=g(f_T(e_i))$ if $i \in T$. 
Every element $m \in M$ is a sum of any coset representative from $m+N$ and some element of $N$, and since $f_S$ and $f_T$ are surjective, $f$ is surjective. 

For any $x \in R_0^{S \sqcup T}$, $x=y+z$ for some $y \in R_0^S$ and $z \in R_0^T$, and $f(x)=f(y+z)=f(y)+f(z)$. If $x=y+z \in \ker(f)$, then $f(y)+f(z)={\bf 0}_M$. Let $h:M \to M/N$ be the quotient map $h(m)=m+N$ for any $m \in M$. then $h(f(y)+f(z))=h(f(y))+h(f(z))=h({\bf 0}_M)$ or $f_T(z)={\bf 0}_{M/N}$. Since $f_T$ is injective, $z={\bf 0}$. Therefore $f(y)={\bf 0}_M$ or $f_S(y)={\bf 0}_M$. Since $f_S$ is injective, $y={\bf 0}$. Thus $x={\bf 0}$ and $f$ is injective.

Since $f$ is bijective, it is an isomorphism, and $M$ is free with ${\rm rank}_I(M)={\rm rank}_I(N)+{\rm rank}_I(M/N).$
\end{proof}

\begin{theorem}
Let $R$ be a commutative ring with identity, let $M$ and $N$ be $R$-modules, and let $f:M \to N$ be an $R$-module homomorphism. If $\ker(f)$ and ${\rm im}(f)$ are free, then $M$ is free and
$${\rm rank}_I(M)={\rm rank}_I(\ker(f))+{\rm rank}_I({\rm im}(f)).$$
\end{theorem}

\begin{proof}
By the First Isomorphism Theorem, $M/\ker(f)$ is isomorphic to ${\rm im}(f)$, and by Theorem \ref{infinite:isoequalrank}, $M/ker(f)$ is free and ${\rm rank}_I(M/\ker(f))={\rm rank}({\rm im}(f))$.
By Theorem \ref{thm:infinite:rank(M/N)}, $${\rm rank}_I(M)={\rm rank}_I(\ker(f))+{\rm rank}_I(M/\ker(f))={\rm rank}_I(\ker(f))+{\rm rank}_I({\rm im}(f)).$$
\end{proof}

\section{Isomorphic Basis of a Module with the Infinite Rank Case}

\begin{definition}
Let $R$ be a commutative ring with identity, let $M$ be an $R$-module, and let $S$ be a set. The set (or multiset) $\{m_i|i \in S\}$ for $m_i \in M$ for all $i \in S$ is an {\bf isomorphic basis} of $M$ if and only if the $R$-module homomorphism $f:R_0^S \to M$ such that $f(e_i)=m_i$ for all $i \in S$ is an isomorphism.
\end{definition}

The standard basis of $R_0^S$, $\{e_i|i \in S\}$, is also an isomorphic basis of $R_0^S$ according to this definition using the identity map $I:R_0^S \to R_0^S$.

Since an $R$-module homomorphism with domain $R_0^S$ is uniquely determined by the outputs of the standard basis, for ${\rm rank}_I(M)=|S|$, there is a one-to-one correspondence between the $R$-module isomorphisms $f:R_0^S \to M$ and the ordered isomorphic bases of $M$ (based on the ordering of $S$).

\begin{theorem}
Let $R$ be a commutative ring with identity, let $M$ be an $R$-module, and let $S$ be a set. Then $M$ is free and ${\rm rank}_I(M)=|S|$ if and only if there exists an isomorphic basis of $M$ containing $|S|$ elements.
\end{theorem}

\begin{proof}
If $M$ is free and ${\rm rank}_I(M)=|S|$, then there exists an $R$-module isomorphism $f:R_0^S \to M$, and $\{f(e_i)|i \in S\}$ is an isomorphic basis of $M$ containing $|S|$ elements.

If $\{m_i|i \in S\}$ is an isomorphic basis of $M$, then the $R$-module homomorphism $f:R_0^S \to M$ defined by $f(e_i)=m_i$ for all $i \in S$ is an isomorphism. Since there exists an isomorphism between $R_0^S$ and $M$, $M$ is free and ${\rm rank}_I(M)=|S|$.
\end{proof}

For free modules $M$, ${\rm rank}_I(M)$ is unique, so we get the following corollary.

\begin{corollary}
Let $R$ be a commutative ring with identity, let $M$ be a free $R$-module. Every isomorphic basis of $M$ has cardinality equal to ${\rm rank}_I(M)$.
\end{corollary}

\begin{definition}
Let $R$ be a commutative ring with identity, let $M$ be an $R$-module, and let $S$ be a set. For the set (or multiset) $\{m_i|i \in S\}$ with $m_i \in M$ for all $i \in S$, an {\bf $R$-linear combination} of $\{m_i|i \in S\}$ is an output of the linear transformation $f:R_0^{S'} \to M$ such that $f(e_i)=m_i$ for all $i \in S'$ for some finite subset $S' \subset S$. 
\end{definition}

This means that an $R$-linear combination is $f(x_1e_{i_1}+\cdots+x_ne_{i_n})=x_1m_{i_1}+\cdots+x_nm_{i_n}$ for some nonnegative integer $n$ and some $x_1,\ldots,x_n \in R$, which is the same as the usual definition.

\begin{theorem}
Let $R$ be a commutative ring with identity, let $M$ be an $R$-module, let $S$ be a set, and let $\mathcal{B}=\{m_i|i \in S\}$ such that $m_i \in V$ for all $i \in S$. Then $\mathcal{B}$ is an isomorphic basis of $M$ if and only if every element in $M$ is a unique $R$-linear combination of $\mathcal{B}$.
\end{theorem}

\begin{proof}
Let $f:R_0^S \to M$ be the $R$-module homomorphism such that $f(e_i)=m_i$ for all $i \in S$.

If $\mathcal{B}$ is an isomorphic basis, then $f$ is an isomorphism and every element in $M$ is a unique output for some input $x_1e_{i_1}+\cdots+x_ke_{i_k} \in R_0^S$, which means it is a unique $R$-linear combination $x_1m_{i_1}+\cdots+x_km_{i_k}$.

If every element in $M$ is a unique $R$-linear combination of $\mathcal{B}$, then $f$ is surjective because every element in $M$ is an output of $f$, and $f$ is injective because $f({\bf 0})={\bf 0}_M$ is unique, which implies $\ker(f)=\{{\bf 0}\}$. Thus $f$ is an isomorphism, and $\mathcal{B}$ is an isomorphic basis.
\end{proof}

\section{Surjective and Injective Sets with the Infinite Rank Case}

\begin{definition}
Let $R$ be a commutative ring with identity, let $M$ be an $R$-module, and let $S$ be a set. For the set (or multiset) $\{m_i|i \in S\}$ with $m_i \in M$ for all $i \in S$, the {\bf $R$-span} of $\{m_i|i \in S\}$ is the image of the $R$-module homomorphism $f:R_0^S \to M$ such that $f(e_i)=m_i$ for all $i \in S$, or in other words, ${\rm span}_{R}(\{m_i|i \in S\})={\rm im}(f)$. 
\end{definition}

This means that ${\rm span}_{R}(\{m_i|i \in S\})$ is the set of all linear combinations of finite subsets of $\{m_i|i \in S\}$, which is the same as the usual definition.

\begin{definition}
Let $R$ be a commutative ring with identity, let $M$ be an $R$-module, and let $S$ be a set. The set (or multiset) $\{m_i|i \in S\}$ with $m_i \in M$ for all $i \in S$ is a {\bf surjective set} in $M$ if and only if the $R$-module homomorphism $f:R_0^S \to M$ such that $f(e_i)=m_i$ for all $i \in S$ is surjective.
\end{definition}

Therefore, $\{m_i|i \in S\}$ is a surjective set in $M$ if and only if ${\rm span}_{R}(\{m_i|i \in S\})=M$. The $R$-module homomorphism $f:R_0^S \to M$ defined by $f(e_i)=m_i$ for all $i \in S$ is also surjective if and only if $\{m_i|i \in S\}$ is a generating set of $M$. Therefore the definitions of surjective set in $M$ and spanning set of $M$ are equivalent. 

\begin{theorem}
Let $R$ be a commutative ring with identity, let $M$ be an $R$-module, let $S$ be a set, and let $\{m_i|i \in S\}$ be a set or multiset with $m_i \in M$ for all $i \in S$. Then $\{m_i|i \in S\}$ is a generating set of $M$ if and only if $\{m_i|i \in S\}$ is a surjective set in $V$.
\end{theorem}

We can prove the results involving spanning sets for surjective sets using the definition of surjective sets.

\begin{theorem} \label{infinite:surjset>=rank(M)}
Let $R$ be a commutative ring with identity, let $M$ be a free $R$-module, let $S$ be a set, and let $\{m_i|i \in S\}$ be a surjective set in $M$. Then $|S| \geq {\rm rank}_I(M)$.
\end{theorem}

\begin{proof}
By definition, the $R$-module homomorphism $f_S:R_0^S \to M$ such that $f_S(e_i)=m_i$ for all $i \in S$ is surjective.
By Theorem \ref{infinite:thm:surj:rank(M)>=rank(N)}, ${\rm rank}_I(R_0^S)=|S| \geq {\rm rank}_I(M)$.
\end{proof}

\begin{definition}
Let $R$ be a commutative ring with identity, let $M$ be an $R$-module, and let $S$ be a set. The set (or multiset) $\{m_i|i \in S\}$ with $m_i \in M$ for all $i \in S$ is an {\bf injective set} in $M$ if and only if the $R$-module homomorphism $f:R_0^S \to M$ such that $f(e_i)=m_i$ for all $i \in S$ is injective.
\end{definition}

If $S=\emptyset$ and $\{m_i|i \in S\}=\emptyset$, then $R_0^S=\{{\bf 0}\}$ and the $R$-module homomorphism $f:\{{\bf 0}\} \to M$ is defined by $f({\bf 0})={\bf 0}_M$, so $f$ is injective. This implies the empty set is an injective set in $M$.

\begin{theorem}
Let $R$ be a commutative ring with identity, let $M$ be an $R$-module, let $S$ be a set, and let $\{m_i|i \in S\}$ be a set or multiset with $m_i \in M$ for all $i \in S$. Then $\{m_i|i \in S\}$ is a linearly independent set in $M$ if and only if $\{m_i|i \in S\}$ is an injective set in $M$.
\end{theorem}

\begin{proof}
Let $f:R_0^S \to M$ be the $R$-module homomorphism defined by $f(e_i)=m_i$ for all $i \in S$. Then $f(x_1e_{i_1}+\cdots+x_ke_{i_k})=x_1m_{i_1}+\cdots+x_km_{i_k}$ for all $x_1e_{i_1}+\cdots+x_ke_{i_k} \in R_0^S$ and $x_1e_{i_1}+\cdots+x_ke_{i_k} \in \ker(f)$ if and only if $x_1m_{i_1}+\cdots+x_km_{i_k}={\bf 0}_M$. 

If $\{m_i|i \in S\}$ is a linearly independent set in $M$, then $x_1m_{i_1}+\cdots+x_km_{i_k}={\bf 0}_M$ implies $x_1e_{i_1}+\cdots+x_ke_{i_k}={\bf 0}$, $\ker(f)=\{{\bf 0}\}$, and $f$ is injective, so $\{m_i|i \in S\}$ is an injective set in $M$. If $\{m_i|i \in S\}$ is an injective set in $M$, then $f$ is injective and $\ker(f)=\{{\bf 0}\}$, and $x_1m_{i_1}+\cdots+x_km_{i_k}={\bf 0}_M$ implies $x_1e_{i_1}+\cdots+x_ke_{i_k}={\bf 0}$.
\end{proof}

\section{Generalized Rank}

\begin{theorem} \label{existssurj}
Let $R$ be a commutative ring with identity, and let $M$ be an $R$-module. There exists a set $S$ and a surjective $R$-module homomorphism $f:R_0^S \to M$.
\end{theorem}

\begin{proof}
Let $S=M$ and define $f:R_0^S \to M$ to be the $R$-module homomorphism such that $f(e_m)=m$ for all $m \in M$.
\end{proof}

\begin{corollary}
Let $R$ be a commutative ring with identity, and let $M$ be an $R$-module. There exists a set $S$ with a minimum cardinality such that there exists a surjective $R$-module homomorphism $f:R_0^S \to M$.
\end{corollary}

\begin{proof}
The set of all sets $S$ such that there exists a surjective $R$-module homomorphism $f:R_0^S \to M$ is nonempty by Theorem \ref{existssurj}. By the Axiom of Choice, the class of cardinal numbers is well-ordered, and there exists a minimum cardinality within this set.
\end{proof}

\begin{theorem}
Let $R$ be a commutative ring with identity, let $M$ be an $R$-module, and let $S$ and $T$ be sets. If $f_S:R_0^S \to M$ is an $R$-module epimorphism and $|S|=|T|$, then there exists an $R$-module epimorphism $f_T:R_0^T \to M$.
\end{theorem}

\begin{proof}
If $|S|=|T|$, then there exists an isomorphism $g:T \to S$. Define the $R$-module homomorphism $f_T:R_0^T \to M$ by $f_T(e_t)=f_S(e_{g(t)})$ for all $t \in T$. Since $g$ is an isomorphism and $f_S$ is an $R$-module epimorphism, then $f_T$ is an $R$-module epimorphism.
\end{proof}

Since every module $M$ is the image of a surjective $R$-module homomorphism with domain $R_0^S$ for some set $S$, we can create the following definition for every $R$-module $M$.

\begin{definition}
Let $R$ be a commutative ring with identity, let $M$ be an $R$-module, and let $S$ be a set. The {\bf generalized rank} of $M$ is $|S|$, ${\rm genrank}(M)=|S|$, if and only if there exists a surjective $R$-module homomorphism $f:R_0^S \to M$ and for any set $T$ such that there exists a surjective $R$-module homomorphism $g:R_0^T \to M$, $|S| \leq |T|$.
\end{definition}

\begin{theorem} \label{thm:rank=genrank}
Let $R$ be a commutative ring with identity, and let $M$ be a free $R$-module. Then ${\rm genrank}(M)={\rm rank}_I(M)$.
\end{theorem}

\begin{proof}
Let $S$ be a set, such that ${\rm rank}_I(M)=|S|$. There exists an $R$-module isomorphism $f:R_0^S \to M$.
Let $T$ be a set, such that there exists a surjective $R$-module homomorphism $g:R_0^T \to M$. 

By Theorem \ref{infinite:thm:surj:rank(M)>=rank(N)}, ${\rm rank}_I(R_0^T)=|T| \geq |S|={\rm rank}_I(M)$. By definition of generalized rank, ${\rm genrank}(M)={\rm rank}_I(M)$.
\end{proof}

\begin{theorem} \label{thm:surj:genrank(M)>=genrank(N)}
Let $R$ be a commutative ring with identity, and let $M$ and $N$ be $R$-modules. If there exists a surjective $R$-module homomorphism $f:M \to N$, then $${\rm genrank}(M) \geq {\rm genrank}(N).$$
\end{theorem}

\begin{proof}
Let $S$ be a set, such that ${\rm genrank}(M)=|S|$. There exists a surjective $R$-module homomorphism $g:R_0^S \to M$. Since $f \circ g:R_0^S \to N$ is surjective, $|S| \geq {\rm genrank}(N)$.
\end{proof}

As a consequence of the previous theorem, we get the following result.

\begin{corollary} \label{cor:iso:genrank(M)=genrank(N)}
Let $R$ be a commutative ring with identity, and let $M$ and $N$ be $R$-modules. If there exists an $R$-module isomorphism $f:M \to N$, then $${\rm genrank}(M) = {\rm genrank}(N).$$
\end{corollary}

\begin{theorem} \label{thm:genrank(M/N)}
Let $R$ be a commutative ring with identity, let $M$ be an $R$-module, and let $N$ be a submodule of $M$. 
$${\rm genrank}(M) \leq {\rm genrank}(N)+{\rm genrank}(M/N).$$
\end{theorem}

\begin{proof}
Let $S$ and $T$ be sets such that ${\rm genrank}(N)=|S|$ and ${\rm genrank}(M/N)=|T|$. There exist surjective $R$-module homomorphisms $f_S:R_0^S \to N$ and $f_T:R_0^T \to M/N$. Define $g:M/N \to M$ to be an injective function that sends each element of $M/N$ to a particular coset representative. Define $f:R_0^{S \sqcup T} \to M$ as the $R$-module homomorphism such that $f(e_i)=f_S(e_i)$ if $i \in S$ and $f(e_i)=g(f_T(e_i))$ if $i \in T$. Every element $m \in M$ is a sum of any coset representative from $m+N$ and some element of $N$, and since $f_S$ and $f_T$ are surjective, $f$ is surjective. 

By Theorem \ref{thm:rank=genrank}, ${\rm genrank}(R_0^{S \sqcup T})=|S \sqcup T|=|S|+|T|$. By Theorem \ref{thm:surj:genrank(M)>=genrank(N)}, $|S|+|T| \geq {\rm genrank}(M)$.
\end{proof}

\begin{theorem}
Let $R$ be a commutative ring with identity, let $M$ and $N$ be $R$-modules, and let $f:M \to N$ be an $R$-module homomorphism.
$${\rm genrank}(M) \leq {\rm genrank}(\ker(f))+{\rm genrank}({\rm im}(f)).$$
\end{theorem}

\begin{proof}
By Theorem \ref{thm:genrank(M/N)}, ${\rm genrank}(M) \leq {\rm genrank}(\ker(f))+{\rm genrank}(M/\ker(f))$. By the First Isomorphism Theorem, $M/\ker(f)$ is isomorphic to ${\rm im}(f)$, and by Corollary \ref{cor:iso:genrank(M)=genrank(N)}, ${\rm genrank}(M/\ker(f))={\rm genrank}({\rm im}(f))$.
\end{proof}

\begin{theorem}
Let $R$ be a commutative ring with identity.
$${\rm genrank}(M) = {\rm genrank}(\ker(f))+{\rm genrank}({\rm im}(f))$$
for any $R$-modules $M$ and $N$ and any $R$-module homomorphism $f:M \to N$ if and only if 
$R$ is a field.
\end{theorem}

\begin{proof}
If $R$ is a field, then ${\rm genrank}={\rm rank}_I=\dim_I$, and $\dim_I(M)=\dim_I(\ker(f))+\dim_I({\rm im}(f))$ by the Rank Nullity Theorem.

For any nonzero $a \in R$, $I=Ra$ is a nonzero ideal and hence nonzero submodule of $R$. Let $f:R \to R/I$ be the quotient map $f(r)=r+I$ for all $r \in R$. As a nonzero submodule, ${\rm genrank}(\ker(f))={\rm genrank}(I)=1$. Also, ${\rm genrank}(R)={\rm rank}_I(R)=1$.

If ${\rm genrank}(R)={\rm genrank}(\ker(f))+{\rm genrank}({\rm im}(f))$, then $${\rm genrank}({\rm im}(f))={\rm genrank}(R/I)=0$$ and $R/I=\{{\bf 0}_{R/I}\}$. Therefore $Ra=R$ and $a$ is invertible. Since this is true for any nonzero $a \in R$, $R$ is a field.
\end{proof}

\end{document}